\documentclass[11pt]{amsart}

\usepackage{amssymb,amsthm,amsmath}
\usepackage[numbers,sort&compress]{natbib}
\usepackage{color}
\usepackage{graphicx}
\usepackage{tikz}
\usepackage[colorlinks,
            linkcolor=blue,
            anchorcolor=blue,
            citecolor=blue
            ]{hyperref}

\let\newpf\proof \let\proof\relax 
\newenvironment{pf}{\newpf[\proofname]}{\qed\endtrivlist}

\newcommand{\ba}{\overline{A}}

\def\be{\begin{equation}}
\def\ee{\end{equation}}

\def\ba{{\begin{align}}}
\def\ea{{\end{align}}}

\def\bm{\begin{matrix}}
\def\em{\end{matrix}}

\def\0{{\mathbf 0}}

\newtheorem{Theorem}{Theorem}[section]
\newtheorem*{Theorem*}{Theorem 1.1}
\newtheorem{Lemma}{Lemma}[section]
\newtheorem{Proposition}{Proposition}[section]
\newtheorem{Corollary}{Corollary}[section]
\newtheorem{Remark}{Remark}[section]

\newtheorem{Definition}{Definition}[section]

\numberwithin{equation}{section}

\theoremstyle{definition}

\renewcommand{\mod}{\operatorname{mod}}

\newcommand{\C}{{\mathbb C}}

\newcommand{\N}{{\mathbb N}}
\newcommand{\Q}{{\mathbb Q}}
\newcommand{\R}{{\mathbb R}}
\newcommand{\T}{{\mathbb T}}

\newcommand{\Z}{{\mathbb Z}}

\newcommand{\la}{\langle}
\newcommand{\ra}{\rangle}

\def\B0{{\bold{0}}}

\catcode`\@=12

\def\Empty{}
\newcommand\oplabel[1]{
  \def\OpArg{#1} \ifx \OpArg\Empty {} \else
    \label{#1}
  \fi}

\newcommand{\comm}[1]{}
\newcommand{\comment}[1]{}

\begin{document}

\title[Reducibility and its applications]{ Reducibility of finitely differentiable quasi-periodic cocycles and its spectral applications }

\author{Ao Cai} \address{
Chern Institute of Mathematics and LPMC, Nankai University, Tianjin 300071, China; and Departmento de Matem\'{a}tica and CMAFCIO, Faculdade de Ci\^{e}ncias, Universidade de Lisboa, Portugal.
}

\email{acai@fc.ul.pt; godcaiao@126.com}

\author{Lingrui Ge} \address{
Department of Mathematics, University of California, Irvine CA, 92717
}
\email{lingruig@uci.edu}
\setcounter{tocdepth}{1}

\begin{abstract}
In this paper, we prove the generic version of Cantor spectrum property for quasi-periodic Schr\"{o}dinger operators with finitely smooth and small potentials, and we also show pure point spectrum for a class of multi-frequency $C^k$ long-range operators on $\ell^2(\Z^d)$. These results are based on reducibility properties of finitely differentiable quasi-periodic $SL(2,\R)$ cocycles. More precisely, we prove that if the base frequency is Diophantine, then a $C^k$ $SL(2,\R)$-valued cocycle is reducible if it is close to a constant cocycle, sufficiently smooth and  the rotation number of it is Diophantine or rational with respect to the frequency.
\end{abstract}

\maketitle
\section{Introduction}
Consider the quasi-periodic Schr\"odinger operator $(\ref{1.1})$ with $C^k$ potential:
\begin{equation}\label{1.1}
(H_{V,\alpha,\theta}x)_n=x_{n+1}+x_{n-1}+V(\theta+n\alpha)x_{n}, n\in \Z,
\end{equation}
where $\theta\in \T^d$ is called the initial phase, $\alpha \in \R^d$ (rationally independent) is called the frequency and $V\in C^k(\T^d,\R)$ is called the potential.
The Aubry dual of $H_{V,\alpha,\theta}$ is called the long-range operator on $\ell^2(\Z^d)$:
\begin{equation}\label{1.2}
(\hat{L}_{V,\alpha,\varphi}u)_m=\sum_{k\in \Z^d}V_k u_{m-k}+2\cos (2\pi(\varphi+\la m,\alpha\ra))u_m, m\in \Z^d.
\end{equation}

In the past forty years, various methods have been developed to study the spectral theory of these two operators due to their close relations to Physics. Among these the most well studied are the structure of the spectrum and the spectral type. Note if the potential is analytic, there has been constant progresses \cite{A0,AJ,AYZ,AYZ2,BG,Eli92,Eli97,GS,J,LYZZ} (consult more references therein), but few was known for the smooth case, the main purpose of this paper is to study some spectral properties of the operator with $C^k$ potential through reducibility methods.


\subsection{Cantor spectrum}
 The geometric structure of spectrum has been one of the central topics for Schr\"odinger operators for several decades. It is conjectured by Simon \cite{Simon} that for generic almost periodic Schr\"odinger operator, the spectrum is a Cantor set. The most intensively studied model is probably the almost Mathieu operator (AMO), the spectrum of which was conjectured - this is the famous ``Ten Martini Problem" - to be a Cantor set. This conjecture has been completely solved by Avila and Jitomirskaya \cite{AJ}. Recently, a stronger version of it, the so-called ``Dry Ten Martini Problem", has been proved for non-critical AMO through reducibility methods by Avila-You-Zhou \cite{AYZ2}. Things become extremely difficult when it comes to Schr\"odinger operators with more general analytic potentials. In the region of positive Lyapunov exponent, Goldstein and Schlag \cite{GS} proved that for any analytic potential, the spectrum is a Cantor set for almost every frequency.  In the perturbative regime, Eliasson \cite{Eli92} proved that for a fixed Diophantine frequency, the spectrum is a Cantor set for generic such potentials in the standard analytic topology $C^{\omega}_h(\T^d, \R)$ with $h>0$. Later, Puig \cite{puig} generalized this result to quasi-periodic Schr\"odinger operators with a fixed Diophantine frequency and potentials in the non-perturbative regime. Note that the topological space given in \cite{Eli92,puig} is the analytic topology. We generalize Eliasson's result to the finitely differentiable topology $C^k(\T^d, \R)$. Before we state our theorem, we give some useful notations first.

For a bounded
analytic (possibly matrix valued) function $F$ defined on $ \{ \theta |  | \Im \theta |< h \}$, let
$
\lvert F\rvert _h=  \sup_{ | \Im \theta |< h } \| F(\theta)\| $ and denote by $C^\omega_{h}(\T^d,*)$ the
set of all these $*$-valued functions ($*$ will usually denote $\R$, $sl(2,\R)$,
$SL(2,\R)$). Also we denote $C^\omega(\T^d,*)=\cup_{h>0}C^\omega_{h}(\T^d,*)$, and set  $C^{k}(\T^{d},*)$  to be the space of $k$ times differentiable with continuous $k$-th derivatives matrix-valued functions. The norms are defined as $$\lVert F \rVert _{k}=\sup_{\substack{
                             \lvert k^{'}\rvert\leqslant k,
                             \theta \in \T^{d}
                          }}\lVert \partial^{k^{'}}F(\theta) \rVert,
$$
and
$$
\lVert F \rVert _{0}=\sup_{\theta \in \T^{d}}\lVert F(\theta) \rVert.
$$

We say $\alpha$ is {\it Diophantine} if there exist $\kappa>0$ and $\tau>d-1$ such that
\begin{equation}\label{dio1}
\alpha \in {\rm DC}_d(\kappa,\tau):=\left\{\alpha \in\R^d:  \inf_{j \in \Z}\left|\la n,\alpha  \ra - j \right|
> \frac{\kappa}{|n|^{\tau}}, \forall \, n\in\Z^d\backslash\{0\} \right\}.
\end{equation}

\begin{Theorem}\label{THM4.2}
Let $\alpha\in {\rm DC}_d(\kappa,\tau)$ and consider $C^k(\T^d, \R)$ with $k\geqslant D_0\tau$ for some large numerical constant $D_0\in \N$, there exists $\epsilon$ depending on $\kappa,\tau,d,k$ such that for a generic potential in
$$
\{V\in C^k(\T^d,\R), \lVert V\rVert_k <\epsilon\}
$$
with respect to the $\lVert\cdot\rVert_k$-topology,
the spectrum of the operator $H_{V,\alpha,\theta}$ is a Cantor set.
\end{Theorem}

\begin{Remark}
A property is $($Baire$)$ generic in $\lVert\cdot\rVert_k$ if the set holding this property contains a residual subset $($a countable intersection of dense open sets$)$ in the $\lVert\cdot\rVert_k$-topology.
\end{Remark}

We review other results on Cantor spectrum for quasi-periodic Schr\"odinger operators with non-analytic potentials. For a class of $C^2$ potentials satisfying a generic transversality condition, Sinai \cite{Sinai} proved that these operators have Cantor spectrum. This was also recently proved by Wang and Zhang \cite{WZ} inspired by L. Young's methods \cite{Young}. For $C^0$ case, Avila-Bochi-Damanik \cite{ABD} showed that for any irrational frequency, and for $C^0$ generic potentials, the spectrum is a Cantor set.

\subsection{Pure point spectrum} The spectral type of quasi-periodic self-adjoint operator is very important in understanding the motion of a particle in a quasi-periodic potential. For the most intensively studied model, the almost Mathieu operator with $V(x)=2\lambda\cos x$, Jitomirskaya \cite{J} proved that if $\lambda>1$ and the frequency is Diophantine, then it has Anderson localization (pure point with exponentially decaying eigenfucntions) for a.e. phase. For general analytic potential with positive Lyapunov exponent, Bourgain and Goldstein \cite{BG} proved that for any fixed phase, it has Anderson localization for a.e. frequency. For more general Gevrey potential, Eliasson \cite{Eli97} proved that for Diophantine frequency, it has pure point for a.e. phase. Recently, based on reducibility method, Avila, You and Zhou \cite{AYZ} solved Jitomirskaya's conjecture \cite{Ji95} on sharp phase transitions of almost Mathieu operator. Note \cite{AYZ} is the measure version of the phase transition conjecture while Jitomirskaya and Liu solved the arithmetic version in \cite{JLiu}. Later, Jitomirskaya and Kachkovskiy \cite{JK} proved that if the Aubry duality of a long-range operator on $\ell^2(\Z^d)$ is an analytic Schr\"odinger operator and the frequency is Diophantine, then it has pure point for a.e. phase. Naturally, people may ask whether the same results hold for the duality of a $C^k$ Schr\"odinger operator? We prove that it is true. Indeed, we can prove a stronger version of it, so-called ``semi-uniformly localized". We first cite the standard definition of ``SULE" introduced in \cite{JL}.

\begin{Definition}\label{semi}\cite{JL}
Assume $H$ is a Schr\"{o}dinger operator acting on $\ell^2(\Z)$ $($this definition also works for $\ell^2(\Z^d)$$)$, it has exponential SULE $($short for ``semi-uniformly localized eigenvectors''$)$ if there exists a constant $\gamma>0$ such that for any $b>0$, there exists a constant $C(b)>0$ such that for any eigenfunction $\psi_s$ one can find $n(s)\in \Z$ such that $\lvert \psi_s(k) \rvert\leqslant C(b)e^{b\lvert n(s)\rvert-\gamma\lvert k-n(s)\rvert }$ for any $k\in \Z$.
\end{Definition}

However, this exponential version of definition does not suit for $C^k$ operators as one can not expect exponential decay of eigenfunction for finitely smooth operators. In $C^k$ case, ``polynomial decay" is much more appropriate:
\begin{Definition}\label{defsemi}
Assume $H$ is a $C^k$ operator acting on $\ell^2(\Z^d)$, it has polynomial SULE if there exists two constants $\tilde{k}>0$ depending on $k$ and $b>0$ independent of $k$ such that for any eigenfunction $u_{s}$, one can find a constant $C>0$ and a site $m(s)\in\Z^d$ such that $|u_{s}(n)|\leq C|m(s)|^b|n-m(s)|^{-\tilde{k}}$ for any $n\in \Z^d$.
\end{Definition}
\begin{Remark}
We can see that ``semi-uniformly localization" contains more information on the eigenfunctions, more precisely on their stability. For more details, one can consult \cite{DJLS}, \cite{JL} and the references therein.
\end{Remark}

With the above definitions, we state our second theorem as follows.
\begin{Theorem}\label{thm4.1}
Let $\alpha\in {\rm DC}_d(\kappa,\tau)$ and $V \in C^k(\T^d, \R)$ with $k\geqslant D_0\tau$ for some large $D_0\in \N$. There exists $\lambda_0$ depending on $\kappa,\tau,d,k,V$ such that if
$
\lvert \lambda \rvert \leqslant \lambda_0,
$
then $\hat{L}_{\lambda V,\alpha,\varphi}$ has pure point spectrum for a.e. $\varphi\in \T$. Furthermore, $\hat{L}_{\lambda V,\alpha,\varphi}$ has polynomial semi-uniformly localized eigenvectors.
\end{Theorem}
\begin{Remark}
This is also an extension of the results in \cite{AYZ} and \cite{JK} in the quasi-periodic case. Moreover, compared to \cite{BG}, the advantage of our method is that we can fix the frequency.
\end{Remark}

It is known that the study of the spectral type of quai-periodic Schr\"odinger operator with lower regularity potential is much more difficult than the analytic case. One of the few results we know is that in the 1980s, Sinai proved Anderson localization for a class of $C^2$ quasi-periodic potentials for any Diophantine frequency in \cite{Sinai}. Later, Fr\"ohlich, Spencer and Wittwer got the same result by different method in \cite{FSW}. For $C^0$ case, Avila and Damanik \cite{AD} proved that for 1-D discrete Schr\"odinger operator with an ergodic continuous potential, there exists a generic set of such potential for which the spectrum of the corresponding operator has no absolutely continuous component.

\subsection{Reducibility}
Note that Schr\"{o}dinger operators \eqref{1.1} are closely related to Schr\"{o}dinger cocycles, where
$$
A(\theta)=S_E^{V}(\theta)=\begin{pmatrix} E-V(\theta) & -1 \\ 1 & 0 \end{pmatrix},
$$
since any solution of the generalized eigenvalue equation $H_{V,\alpha,\theta}u=Eu$ satisfies
$$
A(\theta+n\alpha)\begin{pmatrix} u_{n} \\ u_{n-1} \end{pmatrix}=\begin{pmatrix} u_{n+1} \\ u_{n} \end{pmatrix}.
$$

Therefore, we can analyze the dynamics of  $C^k$ quasi-periodic cocycle $(\alpha,A)\in C^k(\T^d,$ $SL(2,\R))$ and apply it to obtain the spectral properties of the operator \eqref{1.1}.  Reducibility has been widely known as a powerful tool in the study of spectral properties of the spectrum of quasi-periodic Schr\"odinger operators with analytic potentials \cite{AJ,AK,AYZ,AYZ2,Eli92,puig}.  Our method in proving Theorem \ref{THM4.2} and \ref{thm4.1} is a continuation of this powerful method. The philosophy is that nice quantitative reducibility implies nice spectral applications. Here ``quantitative" means that we have very precise estimates on both the conjugation map and the small perturbation.  Before explaining our result precisely, we introduce some basic concepts.

\textbf{Analytic definitions:} Given $h'<h$, we say two cocycles $(\alpha,A_1)$ $(\alpha,A_2)$ $\in \T^d    \times C^{\omega}_h(\T^d,SL(2,\R)): (\theta,v)\mapsto (\theta+\alpha,A_i(\theta)\cdot v), i=1,2$, are $C^{\omega}_{h'}$ conjugated if there exists $Z\in C^{\omega}_{h'}(2\T^d, $ $SL(2,\R))$, such that
$$
Z(\theta+\alpha)A_1(\theta)Z(\theta)^{-1}=A_2(\theta).
$$

An analytic cocycle $(\alpha, A)\in \T^d    \times C^{\omega}_h(\T^d,SL(2,\R))$ is said to be almost reducible if there exist a sequence of conjugations $Z_j\in C^{\omega}_{h_j}(2\T^d, SL(2,\R))$, a sequence of constant matrices $A_j\in SL(2,\R)$ and a sequence of small perturbation $f_j \in C^{\omega}_{h_j}(\T^d, sl(2,\R))$ such that
$$
Z_j(\theta+\alpha)A(\theta)Z_j(\theta)^{-1}=A_j e^{f_j(\theta)}
$$
with
$$
\lvert f_j(\theta)\rvert_{h_j}\rightarrow 0, \ \ j\rightarrow \infty.
$$

Furthermore, we call it weak ($C^{\omega}$) almost reducible if $h_j \rightarrow 0$ and we call it strong ($C^{\omega}_{h,h'}$) almost reducible if $h_j\rightarrow h'>0$. We say $(\alpha, A)$ is $C^{\omega}_{h,h'}$ reducible if there exist a conjugation map $\tilde{Z}\in C^{\omega}_{h'}(2\T^d, $ $SL(2,\R))$ and a constant matrix $\tilde{A} \in SL(2,\R)$ such that
$$
\tilde{Z}(\theta+\alpha)A(\theta)\tilde{Z}(\theta)^{-1}=\tilde{A}(\theta).
$$

\textbf{$C^k$ definitions:} Given $k_1<k$, we say two cocycles $(\alpha,A_1)$ $(\alpha,A_2)\in \T^d    \times C^{k}(\T^d,SL(2,\R)): (\theta,v)\mapsto (\theta+\alpha,A_i(\theta)\cdot v), i=1,2$, are $C^{k_1}$ conjugated if there exists $B\in C^{k_1}(2\T^d, $ $SL(2,\R))$, such that $$
B(\theta+\alpha)A_1(\theta)B(\theta)^{-1}=A_2(\theta).
$$

Cocycle $(\alpha,A)$ is said to be $C^{k,k_1}$ almost reducible, if $A\in C^k(\T^d,SL(2,\R))$ and the $C^{k_1}$-closure of its $C^{k_1}$ conjugacies contains a constant. Moreover, we say $(\alpha,A)$ is  $C^{k,k_1}$ reducible, if $A\in C^k(\T^d,SL(2,\R))$ and its $C^{k_1}$ conjugacies contain a constant.

For quasi-periodic $SL(2,\R)$ cocycle $(\alpha,A)$,  the arithmetic condition on its rotation number $\rho(\alpha, A)$ (see Section 2 $(\ref{rho})$ for basic definition) plays an important role in its reducibility properties. The rotation number $\rho(\alpha, A)$ is said to be {\it Diophantine} with respect to $\alpha$ if there are $\gamma>0$ and $\tau>d-1$ such that $\rho(\alpha, A) \in {\rm DC}_{\alpha}^d(\gamma,\tau)$, where
\begin{equation}\label{dio1}
{\rm DC}_{\alpha}^d(\gamma,\tau):=\left\{\phi \in\R^d:  \inf_{j \in \Z} \left| 2\phi -\la m,\alpha\ra-j\right|
\geqslant \frac{\gamma}{(|m|+1)^{\tau}}, \forall \, m\in\Z^d \right\}.
\end{equation}

The rotation number $\rho(\alpha, A)$ is said to be rational with respect to $\alpha$ if $2\rho(\alpha, A)={\la m_0,\alpha\ra} \mod \Z$ for some $m_0\in \Z^d$.

With the above concepts, our main result is the following:

\begin{Theorem}\label{pro3.1}
Let $A \in SL(2,\R)$, $\alpha\in {\rm DC}_d(\kappa,\tau)$ and $f\in C^k(\T^d, sl(2,\R))$ with $k\geqslant D_0\tau$ for some large $D_0\in \N$. Then there exists $\epsilon$ depending on $\kappa,\tau,d,k,\lVert A\rVert$ such that if $\lVert f(\theta)\rVert_k \leqslant \epsilon$
and $\rho(\alpha, Ae^{f})$ is Diophantine or rational with respect to $\alpha$, then $(\alpha, Ae^{f(\theta)})$ is $C^{k,k_1}$ reducible with $k_1<k$ depending on k.
\end{Theorem}

Reducibility theory has been developed for a long time, we will give a brief review of its history from both local and global side in the analytic case. Note that You-Zhou \cite{youzhou} established a local embedding theorem which builds a bridge between quasi-periodic linear systems and cocycles, showing the equivalence of almost reducibility in the continuous case and the discrete cocycle case. We will review the results according to the  arithmetic condition on the frequency $\alpha$ which is Diophantine or Liouvillean. Recall that we say $\alpha\in \T$ is Liouvillean if $\beta(\alpha)=\limsup \limits_{n\rightarrow \infty}\frac{\ln q_{n+1}}{q_n}>0$ where $q_n$ is the denominator of best rational approximations of $\alpha$.

We first review the reducibility results for local cocycle $(\alpha, Ae^{f(\theta)})$ which is close to constant. Assume that $\alpha\in {\rm DC}_d(\kappa,\tau)$, ``in the perturbative regime'' means $\lvert f(\theta)\rvert_r<\epsilon$ with $\epsilon$ depending on the Diophantine constant $\kappa,\tau$. If $\epsilon$ does not rely on $\kappa,\tau$, it is called ``in the non-perturbative regime''. In the perturbative regime,  Dinaburg and Sinai \cite{DS} first used classical KAM scheme to obtain positive measure (of rotation number) reducibility for continuous quasi-periodic Schr\"odinger equation with small analytic potential. Later, Moser and P\"{o}schel \cite{MP} extended this result to a class of rotation numbers which is rational w.r.t. $\alpha$ by a crucial resonance-cancelation technique. Finally, the breakthrough belongs to Eliasson \cite{Eli92} who proved weak almost reducibility for all energies $E$ and  full measure reducibility for Diophantine frequency and small analytic potential. For strong almost reducibility results, readers may refer to Chavaudret \cite{chavaudret2} and Leguil-You-Zhao-Zhou \cite{LYZZ}. This kind of strong almost reducibility results have many interesting spectral applications, one can consult \cite{LYZZ} and the references therein. In the non-perturbative regime, Puig \cite{puig} used localization method to obtain the non-perturbative version of Eliasson's reducibility theorem. For continuous linear system, Hou-You \cite{houyou} obtained weak almost reducibility for all rotation number, all frequency $\omega=(1,\alpha)\in \T^2$ and small analytic perturbation.

In the Liouvillean case, Avila, Fayad and Krikorian \cite{AFK} obtained that for any $\alpha \in \R \backslash \Q$ and small analytic $f(\theta)$, $(\alpha, Ae^{f(\theta)})$ is full measure rotations reducible (the cocycle can be conjugated to a standard rotation) by ``algebraic conjugacy trick". Meanwhile, Hou-You \cite{houyou} proved full measure rotations reducibility in quasi-periodic linear systems, which can be transferred to the cocycle case by local embedding theorem of You-Zhou \cite{youzhou}. Note the results mentioned above play an essential role in Avila-You-Zhou's solution \cite{AYZ} of Jitomirskaya's conjecture \cite{Ji95}.  More generally, under a non-degeneracy condition, Wang-Zhou \cite{WZhou} used periodic approximation and KAM scheme to get a positive measure diagonalizable result for quasi-periodic $GL(d,\R)$ cocycles which are close to constants.

For the global cocycle $(\alpha, A(\theta))$ with $A\in C^{\omega}(\T,$ $ SL(2,\R))$ and $\alpha \in \T$ (one cannot extend it to $\T^d$ by the technique available now), based on the renormalization scheme, Avila and Krikorian \cite{AK} obtained that for $\alpha$ satisfying some recurrent Diophantine condition, and for almost every $E$, the quasi-periodic Schr\"{o}dinger cocycle is either reducible or non-uniformly hyperbolic (the Lyapunov exponent is positive but the cocycle is not uniformly hyperbolic). Besides, Avila, Fayad and Krikorian \cite{AFK} proved that for irrational $\alpha$,  and  for almost every $E$, the quasi-periodic Schr\"{o}dinger cocycle is either rotations reducible or non-uniformly hyperbolic. Moreover, a recent breakthrough belongs to Avila \cite{A0, A1}, he solved the ``Almost Reducibility Conjecture" i.e. ``subcritical" implies almost reducibility, where ``subcritical" means that there is a uniform subexponential bound on the growth of $\lVert A_n(\theta)\rVert$ ($A_n$ are the transfer matrices) through some band $\lvert \Im \theta \rvert<\delta$.
\subsection{Structure of this paper}In Section 2, some useful definitions and notations are introduced. In Section 3, we prove a quantitative $C^k$ reducibility theorem through KAM, which is a quantitative version of Theorem \ref{pro3.1}. As spectral applications, in Section 4 we use the $C^k$ quantitative reducibility theorem to prove pure point spectrum of the dual long-range operators and generic Cantor spectrum of quasi-periodic Schr\"{o}dinger operators in our settings, which gives Theorem \ref{thm4.1} and Theorem \ref{THM4.2} respectively.
\section{Preliminaries}
\subsection{Hyperbolicity and integrated density of states}
We say the cocycle $(\alpha, A)$ is {\it uniformly hyperbolic} if for every $\theta \in \T^d$, there exists a continuous splitting $\C^2=E^s(\theta)\oplus E^u(\theta)$ such that for some constants $C>0,c>0$, and for every $n\geqslant 0$,
$$
\begin{aligned}
\lvert A_n(\theta)v\rvert \leqslant Ce^{-cn}\lvert v\rvert, & v\in E^s(\theta),\\
\lvert A_{-n}(\theta)v\rvert \leqslant Ce^{-cn}\lvert v\rvert, & v\in E^u(\theta).
\end{aligned}
$$
This splitting is invariant by the dynamics, which means that for every $\theta \in \T^d$, $A(\theta)E^{\ast}(\theta)=E^{\ast}(\theta+\alpha)$, for $\ast=s,u$. The set of uniformly hyperbolic cocycles is open in the $C^0$ topology. In the quasi-periodic Schr\"{o}dinger case, $(\alpha, S_E^V)$ is uniformly hyperbolic if and only if $E\notin \sigma(\alpha,V)$ where $\sigma(\alpha,V)$ denotes the spectrum of $H_{V,\alpha,\theta}$, or equivalently speaking, ``$E$'' lies in some spectral gap \cite{Johnson}.

For the Schr\"{o}dinger operator $H_{V,\alpha,\theta}$, recall that the integrated density of states (IDS) is the function $N:\R \rightarrow [0,1]$ defined by
$$
N(E,V,\alpha)=\int_{\T^{d}}\mu_{V,\alpha,\theta}(-\infty,E]d\theta,
$$
where $\mu_{V,\alpha,\theta}=\mu^{e_{-1}}_{V,\alpha,\theta}+\mu^{e_0}_{V,\alpha,\theta}$ is the spectral measure of $H_{V,\alpha,\theta}$. $\{e_{-1},e_{0}\}$ are the cyclic vectors of $H_{V,\alpha,\theta}$ and $\{e_n\}_{n\in \Z}$ forms a canonical basis of $\ell^2(\Z)$. For more details, readers may refer to Avron-Simon \cite{AS}.

\subsection{Rotation number and degree}
Assume that $A(\cdot): \T^d \rightarrow SL(2,\R))$ is continuous and homotopic to the identity. Then the same is true for the map $F_A:\T^d \times \mathbb{S}^1 \rightarrow \T^d \times \mathbb{S}^1$ with
$$
F_A(x,\omega):=\big(  x+\alpha, \frac{A(x)\cdot \omega}{\lvert A(x)\cdot \omega\rvert}\big).
$$
Thus we can lift $F_A$ to a map $\tilde{F}_A:\T^d\times \R\rightarrow \T^d\times \R$ of the form $\tilde{F}_A(x,y)=(x+\alpha,y+\psi_x(y))$, where for every $x \in \T^d$, $\psi_x(y)$ is $\Z$-periodic. The map $\psi:\T^d \times \R \rightarrow \R$ is called a lift of $A$. Let $\mu$ be any probability measure on $\T^d\times \R$ which is invariant by $\tilde{F}_A$, and whose projection on the first coordinate is given by Lebesgue measure. The number
\begin{equation}\label{rho}
\rho_{(\alpha,A)}:=\int_{\T^d\times \R} \psi_x(y)d\mu(x,y)\mbox{mod}\,\Z
\end{equation}
does not depend on the choices of the lift $\psi$ nor the measure $\mu$. It is called the {\it fibered rotation number} of cocycle $(\alpha,A)$ (readers can consult \cite{JM} for more details).

Let
$$
R_{\phi}:=\begin{pmatrix} \cos2\pi\phi & -\sin2\pi\phi \\ \sin2\pi\phi & \cos2\pi\phi \end{pmatrix},
$$
if $A\in C^0(\T^d,SL(2,\R))$ is homotopic to $\theta\rightarrow R_{\la n,\theta\ra}$ for some $n\in \Z^d$, then we call $n$ the {\it degree} of $A$ and denote it by deg$A$. Moreover, \begin{equation}\label{degree}\deg(AB)=\deg A+\deg B.\end{equation}

Note that the fibered rotation number is invariant under real conjugacies which are homotopic to identity. More generally, if the cocycle $(\alpha,A_1)$ is conjugated to $(\alpha,A_2)$ by $B\in C^0(2\T^d, SL(2,\R))$, i.e. $B(\cdot +\alpha)A_1(\cdot)B^{-1}(\cdot)=A_2(\cdot)$, then
$$
\rho_{(\alpha,A_2)}=\rho_{(\alpha,A_1)}+\frac{\la \deg B,\alpha\ra}{2}.
$$

\subsection{Analytic approximation}
Assume $f \in C^{k}(\T^{d},sl(2,\R))$. By \cite{zehnder}, there exists a sequence $\{f_{j}\}_{j\geqslant 1}$, $f_{j}\in C_{\frac{1}{j}}^{\omega}(\T^{d},sl(2,\R))$ and a universal constant $C^{'}$, such that
\begin{eqnarray} \nonumber \lVert f_{j}-f \rVert_{k} &\rightarrow& 0 , \qquad  j \rightarrow +\infty, \\
\label{2.1}\lvert f_{j}\rvert_{\frac{1}{j}} &\leqslant& C^{'}\lVert f \rVert_{k}, \\   \nonumber \lvert f_{j+1}-f_{j} \rvert_{\frac{1}{j+1}} &\leqslant& C^{'}(\frac{1}{j})^k\lVert f \rVert_{k}.
\end{eqnarray}
Moreover, if $k\leqslant \tilde{k}$ and $f\in C^{\tilde{k}}$, then properties $(\ref{2.1})$ hold with $\tilde{k}$ instead of $k$. That means this sequence is obtained from $f$ regardless of its regularity (since $f_{j}$ is the convolution of $f$ with a map which does not depend on $k$).

\subsection{Aubry duality}
Suppose that the generalized eigenvalue equation of the quasi-periodic Schr\"{o}dinger operator
\begin{equation}
(H_{V,\alpha,\theta}x)_n=x_{n+1}+x_{n-1}+V(\theta+n\alpha)x_{n}=Ex_{n}, n\in \Z
\end{equation}
has a $C^k$ quasi-periodic Bloch wave $x_n=e^{2\pi i n\varphi}\overline{\psi}(\phi+n\alpha)$ for some $\overline{\psi}\in C^k(\T^d,\C)$ and $\varphi \in {[}0,1{)}$. It is easy to see that the Fourier coefficients of $\overline{\psi}(\theta)$ satisfy the following long-range operator:
\begin{equation}
(\hat{L}_{V,\alpha,\varphi}u)_m=\sum_{k\in \Z^d}V_k u_{m-k}+2\cos (2\pi(\varphi+\la m,\alpha\ra))u_m=Eu_m, m\in \Z^d.
\end{equation}
$\hat{L}_{V,\alpha,\varphi}$ is called the dual operator of $H_{V,\alpha,\theta}$.

\section{Full measure reducibility}
In this section, we mainly focus on the reducibility property of the following $C^k$ quasi-periodic $SL(2,\R)$ cocycle:
$$
(\alpha,Ae^{f(\theta)}):\T^{d}\times\R^{2} \rightarrow \T^{d}\times\R^{2};(\theta,v)\mapsto (\theta+\alpha,Ae^{f(\theta)}\cdot v),
$$
where
$f\in C^k(\T^{d},sl(2,\R)), d\in \N^+$ and $\alpha\in {\rm DC}_d(\kappa,\tau)$. However, we are not going to perform KAM method on it directly. Instead, we will analyze the approximating cocycles $\{(\alpha,Ae^{f_{j}(\theta)})\}_{j\geqslant 1}$ first and then obtain the estimates of $(\alpha,Ae^{f(\theta)})$ by analytic approximation $(\ref{2.1})$.

Before we start the proof of Theorem \ref{pro3.1}, we will first cite the Proposition 3.2 of \cite{CCYZ} and give two lemmas to simplify the main proof.
 We first recall some notations given in \cite{CCYZ}. Denote
\begin{equation}
\epsilon_0^{'}(r,r')=\frac{c}{(2\lVert A\rVert)^D}(r-r')^{D\tau},
\end{equation}
and
\begin{equation}\label{mdef}
\epsilon_m=\frac{c}{(2\lVert A \rVert)^D m^{\frac{k}{4}}}, m\in \Z,
\end{equation}
where $c$ depends on $\kappa,\tau,d$ and $D\in \Z$ is a large constant.
For any $k\geqslant 5D\tau$, one can easily compute that

\begin{equation}\label{lowerboundk}\frac{c}{(2\lVert A\rVert)^D m^{\frac{k}{4}}}\leqslant \epsilon_0^{'}(\frac{1}{m},\frac{1}{m^2})
\end{equation}
for any $m\geqslant 10, m\in \Z$.
Denote $l_j=M^{2^{j-1}}$ where both $j\geqslant 1$ and $M>\max\{\frac{(2\lVert A\rVert)^D}{c},10\}$ are integers.

\begin{Proposition}\label{im}
\cite{CCYZ}
Let $\alpha \in {\rm DC}_d(\kappa,\tau)$, $A\in SL(2,\R)$, $\sigma=\frac{1}{10}$, $f\in C^k(\T^{d}$,
$sl(2,\R))$ with $k\geqslant 5D\tau$ and $f_j$ be as in $(\ref{2.1})$, there exists $\epsilon=\epsilon(\kappa,\tau,d,k,\lVert A\rVert)$ such that if $\lVert f\rVert_k\leqslant \epsilon\leqslant \epsilon_0^{'}(\frac{1}{l_{1}},\frac{1}{l_{2}})$, then there exist $B_{l_j}\in C^{\omega}_{\frac{1}{l_{j+1}}}(2\T^{d},SL(2,\R))$, $A_{l_j}\in SL(2,\R)$, $f_{l_j}^{'}\in C^{\omega}_{\frac{1}{l_{j+1}}}(\T^{d},sl(2,\R))$, such that
\begin{equation}\label{fi}
B_{l_j}(\theta+\alpha)(Ae^{f_{l_j}(\theta)})B^{-1}_{l_j}(\theta)=A_{l_j}e^{f_{l_j}^{'}(\theta)},
\end{equation}
with estimates
\begin{equation}\label{estimBlj}
\lvert B_{l_j}(\theta)\rvert_{\frac{1}{l_{j+1}}}\leqslant \epsilon_{l_{j}}^{-\frac{2\sigma}{5}},\ \ \lVert B_{l_j}(\theta)\rVert_0\leqslant \epsilon_{l_{j}}^{-\frac{\sigma}{4}},\ \ \lvert \deg B_{l_j}\rvert\leqslant 4l_{j}\ln\frac{1}{\epsilon_{l_j}},
\end{equation}
\begin{equation}\label{estimflj}
\lvert f_{l_j}^{'}(\theta)\rvert_{\frac{1}{l_{j+1}}}\leqslant \frac{1}{2}\epsilon_{l_{j}}^{\frac{5}{2}},\ \ \lVert A_{l_j}\rVert\leqslant 2\lVert A\rVert.
\end{equation}
Moreover, there exists $\tilde{f}_{l_j}\in C^{k_0}(\T^{d},sl(2,\R))$ with $k_0=[\frac{k}{20}]$ such that
\begin{equation}
B_{l_j}(\theta+\alpha)(Ae^{f(\theta)})B^{-1}_{l_j}(\theta)=A_{l_j}e^{\tilde{f}_{l_j}(\theta)},
\end{equation}
with estimates
\begin{equation}\label{vi}
\lVert \tilde{f}_{l_j}(\theta)\rVert_{k_0} \leqslant \epsilon_{l_j}^2.
\end{equation}
\end{Proposition}

\begin{Remark}
The estimates of $\deg B_{l_j}$ and $(\ref{vi})$ are new, so we will give a brief proof of them.
\end{Remark}
\begin{pf}
Recall that in each KAM step of Proposition 3.1 in \cite{CCYZ}, we have the truncating number $N_j=\frac{2}{\frac{1}{l_j}-\frac{1}{l_{j+1}}}\lvert \ln \epsilon_{l_j}\rvert, j\geqslant 1$. Denote $m_j$ the resonant site of $j$-th KAM step ($m_j=0$ means the $j$-th step is non-resonant), then we have $\lvert m_j\rvert\leqslant N_j$. By the induction process, we have $B_{l_j}=\prod_j^1\tilde{B}_{l_i}$ with $\deg \tilde{B}_{l_i} = m_i$, thus by $(\ref{degree})$, $\deg B_{l_j}=\sum_{i=1}^{j}m_i,$ and
$$
\lvert \deg B_{l_j}\rvert \leqslant \sum_{i=1}^{j}\lvert m_i\rvert \leqslant \sum_{i=1}^{j}N_j\leqslant \sum_{i=1}^{j}\frac{2}{\frac{1}{l_j}-\frac{1}{l_{j}^2}}\lvert \ln \epsilon_{l_j}\rvert \leqslant 4l_j \ln \frac{1}{\epsilon_{l_j}},
$$
which gives $(\ref{estimBlj})$.

By $(\ref{fi})$, we have
$$
B_{l_j}(\theta+\alpha)(Ae^{f(\theta)})B^{-1}_{l_j}(\theta)=A_{l_j}e^{f_{l_j}^{'}(\theta)}+B_{l_j}(\theta+\alpha)(Ae^{f(\theta)}-Ae^{f_{l_j}(\theta)})B^{-1}_{l_j}(\theta).
$$

Denote $A_{l_j}e^{f_{l_j}^{'}(\theta)}+B_{l_j}(\theta+\alpha)(Ae^{f(\theta)}-Ae^{f_{l_j}(\theta)})B^{-1}_{l_j}(\theta)=A_{l_j}e^{\tilde{f}_{l_j}(\theta)},$ then
$$
\lVert \tilde{f}_{l_j}(\theta)\rVert_{k_0}\leqslant \lVert f_{l_j}^{'}(\theta)\rVert_{k_0}+\lVert A_{l_j}^{-1}\rVert \lVert B_{l_j}(\theta+\alpha)(Ae^{f(\theta)}-Ae^{f_{l_j}(\theta)})B^{-1}_{l_j}(\theta)\rVert_{k_0}.
$$

By the same argument of Theorem 3.1 in \cite{CCYZ}, by Cauchy estimates, for $k_0 \in \N$ and $k_0\leqslant k $, we have
\begin{flalign*}
\lVert f_j-f_{j+1}\rVert_{k_0}& \leqslant \sup_{\substack{
                             l\leqslant k_0,
                             \theta \in \T^{d}
                          }}\lVert (\partial_{\theta_1}^{l_1}+\cdots+\partial_{\theta_d}^{l_d})(f_j(\theta)-f_{j+1}(\theta))\rVert \\
                          &\leqslant (k_0)!(j+1)^{k_0}\lvert f_j-f_{j+1}\rvert_{\frac{1}{j+1}}\\
                          &\leqslant (k_0)!(j+1)^{k_0}\times \frac{c}{(2\lVert A\rVert)^D j^{k}}\\
                          &\leqslant \frac{C_1}{j^{k-k_0}},
\end{flalign*}
where $C_1$ is independent of $j$.

By a simple integration we get
$$
\lVert f(\theta)-f_{l_j}(\theta)\rVert_{k_0}\leqslant \frac{C_1}{l_j^{k-k_0-1}}.
$$

Similarly by Cauchy estimates, we have
\begin{flalign*}
\lVert f_{l_j}^{'}(\theta)\rVert_{k_0}& \leqslant (k_0)!(l_{j+1})^{k_0}\lvert f_{l_j}^{'}(\theta)\rvert_{\frac{1}{l_{j+1}}}\\
&\leqslant (k_0)!(l_{j})^{2k_0}\times \frac{1}{2}\epsilon_{l_j}^{\frac{5}{2}}\\
&\leqslant (k_0)!(l_{j})^{2k_0}\times \frac{1}{2}\times (\frac{c}{(2\lVert A\rVert)^D {l_j}^{\frac{k}{4}}})^{\frac{5}{2}}\\
&\leqslant \frac{C_2}{l_j^{\frac{5k}{8}-2k_0}},
\end{flalign*}
where $C_2$ is independent of $j$.
\begin{flalign*}
\lVert B_{l_j}(\theta)\rVert_{k_0}&\leqslant (k_0)!(l_{j+1})^{k_0}\lvert B_{l_j}(\theta)\rvert_{\frac{1}{l_{j+1}}}\\
&\leqslant (k_0)!(l_{j})^{2k_0}\times \epsilon_{l_{j}}^{-\frac{2\sigma}{5}}\\
&\leqslant (k_0)!(l_{j})^{2k_0}\times (\frac{c}{(2\lVert A\rVert)^D {l_j}^{\frac{k}{4}}})^{-\frac{2\sigma}{5}}\\
&\leqslant C_3 l_j^{\frac{k}{100}+2k_0},
\end{flalign*}
where $C_3$ is independent of $j$.
Thus we have
\begin{flalign*}
\lVert \tilde{f}_{l_j}(\theta)\rVert_{k_0} &\leqslant \lVert f_{l_j}^{'}(\theta)\rVert_{k_0}+\lVert A_{l_j}^{-1}\rVert \lVert B_{l_j}(\theta+\alpha)(Ae^{f(\theta)}-Ae^{f_{l_j}(\theta)})B^{-1}_{l_j}(\theta)\rVert_{k_0}.\\
&\leqslant \frac{C_2}{l_j^{\frac{5k}{8}-2k_0}}+C_4 l_j^{\frac{k}{50}+4k_0}\times l_j^{-k+k_0+1}\\
&\leqslant \frac{C_2}{l_j^{\frac{5k}{8}-2k_0}}+\frac{C_4}{l_j^{\frac{49k}{50}-5k_0-1}},
\end{flalign*}
where $C_4$ is independent of $j$.

Thus if we pick $k_0\leqslant \frac{k}{20}$, then
$$
\lVert \tilde{f}_{l_j}(\theta)\rVert_{k_0} \leqslant \epsilon_{l_j}^2,
$$
which gives $(\ref{vi})$.
\end{pf}

Now, we will give our first Lemma concerning positive measure reducibility, i.e. for fixed $\gamma>0,\tau>d-1$, the cocycle $(\alpha,Ae^{f(\theta)})$ is reducible provided $f(\theta)$ is sufficiently small (depending on $\gamma,\tau$).
\begin{Lemma}\label{lem3.1}
Let $\alpha\in {\rm DC}_d(\kappa,\tau)$, $A \in SL(2,\R)$ and $f\in C^k(\T^d, sl(2,\R))$ with $k\geqslant 5D\tau$. Assume that $\rho(\alpha,Ae^f) \in {\rm DC}_{\alpha}^d(\gamma,\tau)$, there exists $T=T(\tau)$ and $\epsilon_1=\epsilon_1(\gamma,\kappa,\tau,d,k,\lVert A\rVert)$ such that if
\begin{equation}\label{sm1}
\lVert f\rVert_k\leqslant \epsilon_1\leqslant T(\tau)\gamma^{11}\epsilon_0^{'}(\frac{1}{l_{1}},\frac{1}{l_{2}}),
\end{equation}
then there exists $B_1\in C^{\tilde{k}}(\T^d, SL(2,\R))$ with $\tilde{k}= [\frac{k}{20}]$ such that
\begin{equation}\label{main11}
B_1(\theta+\alpha)(Ae^{f(\theta)})B_1^{-1}(\theta)=A_1,
\end{equation}
with
\begin{equation}\label{main12}
\lVert B_1(\theta)\rVert_{\tilde{k}}< 1+\delta,
\end{equation}
where $A_1\in SL(2,\R)$ is a constant matrix and $\delta$ is a small constant depending on $\gamma,\kappa,\tau,d,k,
\lVert A\rVert$.
\end{Lemma}
\begin{pf}
We prove this by induction. Let $\{f_{l_j}\}_{j\geqslant 1}$ be the approximating subsequence as above.

\textbf{First step:} By $(\ref{sm1})$ and Proposition \ref{im}, we have
$$
B_{l_1}(\theta+\alpha)(Ae^{f_{l_1}(\theta)})B^{-1}_{l_1}(\theta)=A_{l_1}e^{f_{l_1}^{'}(\theta)}.
$$

Moreover, since $\epsilon_{l_1}\leqslant T(\tau)\gamma^{11}\epsilon_0^{'}(\frac{1}{l_{1}},\frac{1}{l_{2}})$ and $\rho(\alpha,Ae^f) \in {\rm DC}_{\alpha}^d(\gamma,\tau)$, for $\lvert m\rvert \leqslant \frac{2}{\frac{1}{l_{1}}-\frac{1}{l_{2}}}\ln\frac{1}{\epsilon_{l_{1}}}$, we have
\begin{flalign*}
&\lVert 2\rho(\alpha, A)-\la m,\alpha\ra\rVert_{\R/\Z}\\
\geqslant &\lVert 2\rho(\alpha, Ae^{f})-\la m,\alpha\ra\ \rVert_{\R/\Z} -\lvert 2\rho(\alpha, Ae^{f})-2\rho(\alpha, A)\rvert\\
\geqslant &\frac{\gamma}{(\lvert m\rvert+1)^{\tau}}-2\epsilon_{l_{1}}^{\frac{1}{2}}\\
\geqslant &\frac{\gamma}{(\frac{2}{\frac{1}{l_{1}}-\frac{1}{l_{2}}}\ln\frac{1}{\epsilon_{l_{1}}}+1)^{\tau}}-2\epsilon_{l_{1}}^{\frac{1}{2}}\\
\geqslant & 2\epsilon_{l_{1}}^{\frac{1}{10}}-2\epsilon_{l_{1}}^{\frac{1}{2}}\\
> & \epsilon_{l_{1}}^{\frac{1}{10}},
\end{flalign*}
which means the first step is non-resonant with
$$
B_{l_{1}}=\tilde{B}_{l_{1}}\circ B_{l_{0}}, \, B_{l_{0}}:= Id,\ \
\lvert \tilde{B}_{l_1}(\theta)-Id\rvert_{\frac{1}{l_2}}\leqslant \epsilon_{l_1}^{\frac{1}{2}}
$$
and thus $\deg B_{l_{1}}=0$.

\textbf{Induction step:} assume that for $j\leqslant j_0$, we have
$$
B_{l_j}(\theta+\alpha)(Ae^{f_{l_j}(\theta)})B^{-1}_{l_j}(\theta)=A_{l_j}e^{f_{l_j}^{'}(\theta)},
$$
which is equivalent to
\begin{equation}\label{rot}
B_{l_j}(\theta+\alpha)(Ae^{f(\theta)})B^{-1}_{l_j}(\theta)=A_{l_j}e^{f_{l_j}^{'}(\theta)}+B_{l_j}(\theta+\alpha)(Ae^{f(\theta)}-Ae^{f_{l_j}(\theta)})B^{-1}_{l_j}(\theta),
\end{equation}
with estimates
\begin{equation}\label{b}
B_{l_{j}}=\tilde{B}_{l_{j}}\circ B_{l_{j-1}},\ \
\lvert \tilde{B}_{l_j}(\theta)-Id\rvert_{\frac{1}{l_{j+1}}}\leqslant \epsilon_{l_j}^{\frac{1}{2}}, \ \ \deg B_{l_{j-1}}=0.
\end{equation}

Notice that $(\ref{estimBlj})$, $(\ref{estimflj})$ and $(\ref{rot})$ gives
$$
\lvert \rho(\alpha,Ae^f)+\frac{\la \mbox{deg}B_{l_j},\alpha\ra}{2} -\rho(\alpha,A_{l_j})\rvert \leqslant \epsilon_{l_{j+1}}^{\frac{1}{2}},
$$
and $(\ref{b})$ implies $\mbox{deg}B_{l_j}=0$. Thus we have
\begin{equation}\label{rott2}
\lvert \rho(\alpha,Ae^f)-\rho(\alpha,A_{l_j})\rvert \leqslant \epsilon_{l_{j+1}}^{\frac{1}{2}}.
\end{equation}

Consequently, by $(\ref{lowerboundk})$, $(\ref{rott2})$ and the Diophantine condition on $\rho(\alpha, Ae^{f})$, for $\lvert m\rvert \leqslant \frac{2}{\frac{1}{l_{j+1}}-\frac{1}{l_{j+2}}}\ln\frac{1}{\epsilon_{l_{j+1}}}$, we have
\begin{flalign*}
&\lVert 2\rho(\alpha,A_{l_j})-\la m,\alpha\ra\rVert_{\R/\Z}\\
\geqslant & \lVert 2\rho(\alpha,Ae^f)-\la m,\alpha\ra \rVert_{\R/\Z}- \lvert 2\rho(\alpha,Ae^f)-2\rho(\alpha,A_{l_j})\rvert \\
\geqslant & \frac{\gamma}{(\lvert m\rvert+1)^{\tau}}-2\epsilon_{l_{j+1}}^{\frac{1}{2}}\\
\geqslant & \frac{\gamma}{(\frac{2}{\frac{1}{l_{j+1}}-\frac{1}{l_{j+2}}}\ln\frac{1}{\epsilon_{l_{j+1}}}+1)^{\tau}}-2\epsilon_{l_{j+1}}^{\frac{1}{2}}\\
\geqslant & 2\epsilon_{l_{j+1}}^{\frac{1}{10}}-2\epsilon_{l_{j+1}}^{\frac{1}{2}}\\
\geqslant & \epsilon_{l_{j+1}}^{\frac{1}{10}}, \ \  \forall \,j\leqslant j_0.
\end{flalign*}

This means the $(j_0+1)$-th step is also non-resonant with estimates
$$
B_{l_{j_0+1}}=\tilde{B}_{l_{j_0+1}}\circ B_{l_{j_0}},\ \
\lvert \tilde{B}_{l_{j_0+1}}(\theta)-Id\rvert_{\frac{1}{l_{j_0+2}}}\leqslant \epsilon_{l_{j_0+1}}^{\frac{1}{2}}.
$$

To conclude, we have $\forall \, j\geqslant 1$,
$$
B_{l_j}(\theta+\alpha)(Ae^{f_{l_j}(\theta)})B^{-1}_{l_j}(\theta)=A_{l_j}e^{f_{l_j}^{'}(\theta)},
$$
with
$$
B_{l_{j}}=\tilde{B}_{l_{j}}\circ B_{l_{j-1}},\ \
\lvert \tilde{B}_{l_j}(\theta)-Id\rvert_{\frac{1}{l_{j+1}}}\leqslant \epsilon_{l_j}^{\frac{1}{2}}.
$$

Denote $B_1=\lim_{j\rightarrow \infty}B_{l_j}, A_1=\lim_{j\rightarrow \infty}A_{l_j}$, then $(\ref{main11})$ holds. By Cauchy estimates, for $\tilde{k} \in \N$ and $\tilde{k}\leqslant k $, we have
\begin{flalign*}
\lVert \tilde{B}_{l_{j}}(\theta)-Id\rVert_{\tilde{k}}& \leqslant \sup_{\substack{
                             \lvert l\rvert \leqslant \tilde{k},
                             \theta \in \T^{d}
                          }}\lVert (\partial_{\theta_1}^{l_1}+\cdots+\partial_{\theta_d}^{l_d})(\tilde{B}_{l_{j}}(\theta)-Id)\rVert \\
                          &\leqslant (\tilde{k})!(l_{j+1})^{\tilde{k}}\lvert \tilde{B}_{l_{j}}-Id\rvert_{\frac{1}{l_{j+1}}}\\
                          &\leqslant (\tilde{k})!(l_{j})^{2\tilde{k}}\times \frac{C}{{l_{j}}^{\frac{k}{8}}}\\
                          &\leqslant \frac{C_1}{l_{j}^{\frac{k}{8}-2\tilde{k}}},
\end{flalign*}
where $C_1$ does not depend on $j$. Thus if we pick $\tilde{k}\leqslant \frac{k}{20}$, we have
$$
\lVert \tilde{B}_{l_{j}}-Id\rVert_{\tilde{k}}\leqslant \frac{C_1}{l_{j}^{\frac{k}{40}}},
$$
and then
$$
\lVert \tilde{B}_{l_{j}}\rVert_{\tilde{k}}\leqslant 1+\frac{C_1}{l_{j}^{\frac{k}{40}}}.
$$

Therefore, since we pick $l_1=M$ sufficiently large, then
$$
\lVert B_1\rVert_{\tilde{k}}=\lVert \prod_{j}\tilde{B}_{l_{j}}\rVert_{\tilde{k}}\leqslant \prod_{j}(1+\frac{C_1}{l_{j}^{\frac{k}{40}}})\leqslant 1+\frac{2C_1}{M^{\frac{k}{40}}}< 1+\delta.
$$
thus $(\ref{main12})$ holds. This finishes the proof of Lemma $\ref{lem3.1}$.
\end{pf}

Note that Lemma $\ref{lem3.1}$ is used to prove Theorem $\ref{pro3.1}$ when the rotation number is Diophantine with respect to $\alpha$ while Lemma $\ref{lem3.2}$ below is useful to the proof when the rotation number is rational with respect to $\alpha$.
\begin{Lemma}\label{lem3.2}
Let $\alpha\in {\rm DC}_d(\kappa,\tau)$, $A \in SL(2,\R)$ and $f\in C^k(\T^d, sl(2,\R))$ with $k\geqslant 5D\tau$. Assume that $\rho(\alpha,Ae^f)=0$, there exists $T=T(\tau)$ and $\epsilon_2=\epsilon_2(\kappa,\tau,d,k,
\lVert A\rVert)$ such that if
\begin{equation}\label{sm2}
\lVert f\rVert_k\leqslant \epsilon_2\leqslant T(\tau)\kappa^{11}\epsilon_0^{'}(\frac{1}{l_{1}},\frac{1}{l_{2}}),
\end{equation}
then there exists $B_2\in C^{\tilde{k}}(\T^d, SL(2,\R))$ with $\tilde{k}= [\frac{k}{20}]$ such that
\begin{equation}\label{main3}
B_2(\theta+\alpha)(Ae^{f(\theta)})B_2^{-1}(\theta)=A_2,
\end{equation}
with
\begin{equation}\label{main4}
\lVert B_2(\theta)\rVert_{\tilde{k}}< 1+\delta,
\end{equation}
where $A_2\in SL(2,\R)$ is a constant matrix and $\delta$ is a small constant depending on $\kappa,\tau,d,k,
\lVert A\rVert$.
\end{Lemma}
\begin{pf}
\textbf{First step:} By $(\ref{sm2})$ and Proposition \ref{im}, we have
\begin{equation}\label{ini}
B_{l_1}(\theta+\alpha)(Ae^{f_{l_1}(\theta)})B^{-1}_{l_1}(\theta)=A_{l_1}e^{f_{l_1}^{'}(\theta)}.
\end{equation}

Moreover, since $\epsilon_{l_1}\leqslant T(\tau)\kappa^{11}\epsilon_0^{'}(\frac{1}{l_{1}},\frac{1}{l_{2}})$ and $\rho(\alpha,Ae^f)=0$, thus for $\lvert m\rvert \leqslant \frac{2}{\frac{1}{l_{1}}-\frac{1}{l_{2}}}\ln\frac{1}{\epsilon_{l_{1}}}$, we have
\begin{flalign*}
&\lVert 2\rho(\alpha, A)-\la m,\alpha\ra\rVert_{\R\backslash\Z}\\
\geqslant &\lVert 2\rho(\alpha, Ae^{f})-\la m,\alpha\ra \rVert_{\R/\Z} -\lvert 2\rho(\alpha, Ae^{f})-2\rho(\alpha, A)\rvert\\
\geqslant &\lVert \la m,\alpha\ra \rVert_{\R/\Z} -\lvert 2\rho(\alpha, Ae^{f})-2\rho(\alpha, A)\rvert\\
\geqslant &\frac{\kappa}{\lvert m\rvert^{\tau}}-2\epsilon_{l_{1}}^{\frac{1}{2}}\\
\geqslant &\frac{\kappa}{(\frac{2}{\frac{1}{l_{1}}-\frac{1}{l_{2}}}\ln\frac{1}{\epsilon_{l_{1}}})^{\tau}}-2\epsilon_{l_{1}}^{\frac{1}{2}}\\
\geqslant & 2\epsilon_{l_{1}}^{\frac{1}{10}}-2\epsilon_{l_{1}}^{\frac{1}{2}}\\
> & \epsilon_{l_{1}}^{\frac{1}{10}},
\end{flalign*}
which means the first step is non-resonant with
$$
B_{l_{1}}=\tilde{B}_{l_{1}}\circ B_{l_{0}}, \, B_{l_{0}}:= Id,\ \
\lvert \tilde{B}_{l_1}(\theta)-Id\rvert_{\frac{1}{l_2}}\leqslant \epsilon_{l_1}^{\frac{1}{2}},
$$
and thus $\deg B_{l_{1}}=0$.

\textbf{Induction step:} assume that for $j\leqslant j_0$, we have
$$
B_{l_j}(\theta+\alpha)(Ae^{f_{l_j}(\theta)})B^{-1}_{l_j}(\theta)=A_{l_j}e^{f_{l_j}^{'}(\theta)},
$$
which is equivalent to
\begin{equation}\label{rot2}
B_{l_j}(\theta+\alpha)(Ae^{f(\theta)})B^{-1}_{l_j}(\theta)=A_{l_j}e^{f_{l_j}^{'}(\theta)}+B_{l_j}(\theta+\alpha)(Ae^{f(\theta)}-Ae^{f_{l_j}(\theta)})B^{-1}_{l_j}(\theta),
\end{equation}
with estimates
\begin{equation}\label{b2}
B_{l_{j}}=\tilde{B}_{l_{j}}\circ B_{l_{j-1}},\ \
\lvert \tilde{B}_{l_j}(\theta)-Id\rvert_{\frac{1}{l_{j+1}}}\leqslant \epsilon_{l_j}^{\frac{1}{2}}, \ \ \deg B_{l_{j-1}}=0.
\end{equation}

Notice that $(\ref{estimBlj})$, $(\ref{estimflj})$ and $(\ref{rot2})$ gives
$$
\lvert \rho(\alpha,Ae^f)+\frac{\la \mbox{deg}B_{l_j},\alpha\ra}{2} -\rho(\alpha,A_{l_j})\rvert \leqslant \epsilon_{l_{j+1}}^{\frac{1}{2}},
$$
and $(\ref{b2})$ implies $\mbox{deg}B_{l_j}=0$. Thus we have
\begin{equation}\label{rott3}
\lvert \rho(\alpha,Ae^f)-\rho(\alpha,A_{l_j})\rvert \leqslant \epsilon_{l_{j+1}}^{\frac{1}{2}}.
\end{equation}

By $(\ref{lowerboundk})$, $(\ref{rott3})$ and $\rho(\alpha, Ae^{f})=0$, for $\lvert m\rvert \leqslant \frac{2}{\frac{1}{l_{j+1}}-\frac{1}{l_{j+2}}}\ln\frac{1}{\epsilon_{l_{j+1}}}$, we have
\begin{flalign*}
&\lVert 2\rho(\alpha, A_{l_j})-\la m,\alpha\ra\rVert_{\R/\Z}\\
\geqslant &\lVert 2\rho(\alpha, Ae^{f})-\la m,\alpha\ra \rVert_{\R/\Z} -\lvert 2\rho(\alpha, Ae^{f})-2\rho(\alpha, A_{l_j})\rvert\\
\geqslant &\lVert \la m,\alpha\ra \rVert_{\R/\Z} -\lvert 2\rho(\alpha, Ae^{f})-2\rho(\alpha, A_{l_j})\rvert\\
\geqslant &\frac{\kappa}{\lvert m\rvert^{\tau}}-2\epsilon_{l_{j+1}}^{\frac{1}{2}}\\
\geqslant &\frac{\kappa}{(\frac{2}{\frac{1}{l_{j+1}}-\frac{1}{l_{j+2}}}\ln\frac{1}{\epsilon_{l_{j+1}}})^{\tau}}-2\epsilon_{l_{j+1}}^{\frac{1}{2}}\\
\geqslant & 2\epsilon_{l_{j+1}}^{\frac{1}{10}}-2\epsilon_{l_{j+1}}^{\frac{1}{2}}\\
> & \epsilon_{l_{j+1}}^{\frac{1}{10}}, \ \ \forall \,j\leqslant j_0,
\end{flalign*}
which means the $(j_0+1)$-th step is non-resonant. By the same process of Lemma \ref{lem3.1}, we obtain $(\ref{main3})$ and $(\ref{main4})$. This finishes our proof.
\end{pf}

With the above Proposition and Lemmas in hand, we are ready to prove the following quantitative version of Theorem \ref{pro3.1}.

\begin{Theorem}\label{quanthm}
Let $\alpha\in {\rm DC}_d(\kappa,\tau)$, $A \in SL(2,\R)$ and $f\in C^k(\T^d, sl(2,\R))$ with $k\geqslant D_0\tau$ for some large $D_0\in \N$. Then there exists $\epsilon$ depending on $\kappa,\tau,d,k,\lVert A\rVert$ such that if
\begin{equation}\label{sm3}
\lVert f(\theta)\rVert_k \leqslant \epsilon\leqslant \epsilon_0^{'}(\frac{1}{l_{1}},\frac{1}{l_{2}}),
\end{equation}
and
\begin{itemize}
\item if $\rho(\alpha, Ae^f)$ is Diophantine with respect to $\alpha:$ $\rho(\alpha, Ae^f)\in {\rm DC}_{\alpha}^d(\gamma,\tau)$, then there exists two constants $C_1=\gamma^{-20}C(\kappa,\tau,
d,k,\lVert A\rVert)$, $C_2=\gamma^{-\frac{100}{k}}C(\kappa,\tau,d,k,\lVert A\rVert)$ and $\tilde{B}_1 \in C^{k_1}$ $(2\T^d,SL(2,\R))$ with $k_1=[\frac{k}{400}]$ such that
\begin{equation}\label{esta}
\tilde{B}_1^{-1}(\theta+\alpha)Ae^{f(\theta)}\tilde{B}_1(\theta)=R_{\phi}\in SL(2,\R), \, \phi\notin \Z,
\end{equation}
with estimates:
\begin{equation}\label{estb}
\lVert \tilde{B}_1\rVert_{k_1}\leqslant C_1, \ \ \lvert \deg \tilde{B}_1\rvert \leqslant C_2.
\end{equation}
\item if $\rho(\alpha, Ae^f)$ is rational with respect to $\alpha:$ $2\rho(\alpha, A)={\la m_0,\alpha\ra} \mod \Z$ for some $m_0\in \Z^d$, then there exists $\tilde{B}_2 \in C^{k_1}(2\T^d,SL(2,\R))$ with $k_1=[\frac{k}{400}]$ such that
\begin{equation}\label{estbb}
\tilde{B}_2^{-1}(\theta+\alpha)Ae^{f(\theta)}\tilde{B}_2(\theta)=\bar{A}\in SL(2,\R), \,  \rho(\alpha,\bar{A})=0.
\end{equation}
\end{itemize}
\end{Theorem}
\begin{pf}
\textbf{(Diophantine case)} By $(\ref{sm3})$ and Proposition \ref{im}, there exists $B_{l_j}\in C^{\omega}_{\frac{1}{l_{j+1}}}(2\T^{d},SL(2,\R))$, $A_{l_j}\in SL(2,\R)$ and $\tilde{f}_{l_j}\in C^{k_0}(\T^{d},sl(2,\R))$ with $k_0= [\frac{k}{20}]$ such that
\begin{equation}
B_{l_j}(\theta+\alpha)(Ae^{f(\theta)})B^{-1}_{l_j}(\theta)=A_{l_j}e^{\tilde{f}_{l_j}(\theta)},
\end{equation}
with estimates
\begin{equation}\label{ti}
\lvert B_{l_j}(\theta)\rvert_{\frac{1}{l_{j+1}}}\leqslant \epsilon_{l_{j}}^{-\frac{2\sigma}{5}},\, \lvert \deg B_{l_j}\rvert\leqslant 4l_{j}\ln\frac{1}{\epsilon_{l_j}},\, \lVert \tilde{f}_{l_j}(\theta)\rVert_{k_0} \leqslant \epsilon_{l_j}^2.
\end{equation}

Since $\rho(\alpha, Ae^f)\in {\rm DC}_{\alpha}^d(\gamma,\tau)$, for any $m\in \Z^d$, we have
\begin{flalign*}
& \lVert 2\rho(\alpha, A_{l_j}e^{\tilde{f}_{l_j}(\theta)}) - \la m,\alpha\ra\rVert_{\R/\Z}\\
= & \lVert 2\rho(\alpha, Ae^f)+\la \deg B_{l_j},\alpha\ra -\la m,\alpha\ra\rVert_{\R/\Z}\\
\geqslant & \frac{\gamma}{(\lvert m- \deg B_{l_j}\rvert+1)^{\tau}}\\
\geqslant & \frac{\gamma (1+\lvert \deg B_{l_j}\rvert )^{-\tau}}{(\lvert m\rvert+1)^{\tau}},
\end{flalign*}
which implies $\rho(\alpha, A_{l_j}e^{\tilde{f}_{l_j}(\theta)})\in {\rm DC}_{\alpha}^d(\gamma(1+\lvert \deg B_{l_j}\rvert )^{-\tau},\tau)$. Let $j$ be the smallest integer such that
\begin{equation}\label{pickj}
\epsilon_{l_j}^2 \leqslant T(\tau)(\frac{\gamma}{(1+4l_{j}\ln\frac{1}{\epsilon_{l_j}} )^{\tau}})^{11}\epsilon_0^{'}(\frac{1}{l_{1}},\frac{1}{l_{2}}),
\end{equation}
where $T(\tau)$ is defined in Lemma \ref{lem3.1} (such $j=j(\gamma)$ exists since we choose $k\geqslant D_0\tau$ with $D_0$ large enough). Thus by Lemma \ref{lem3.1}, there exists $B_1\in C^{k_1}(\T^d, SL(2,\R))$ with $k_1= [\frac{k_0}{20}]$ such that
\begin{equation}\label{main1}
B_1(\theta+\alpha)(A_{l_j}e^{\tilde{f}_{l_j}(\theta)})B_1^{-1}(\theta)=A_1,
\end{equation}
where $A_1\in SL(2,\R)$ is a constant matrix and $B_1(\theta)$ is close to identity (thus $\deg B_1=0$). Denote
$$
\tilde{B_1}=B_1\circ B_{l_j},
$$
by $(\ref{ti})$ and Cauchy estimates, we have
\begin{equation}\label{cau}
\lVert B_{l_j}\rVert_{k_1}\leqslant (k_1)!(l_{j+1})^{k_1}\lvert B_{l_j}\rvert_{\frac{1}{l_{j+1}}}\leqslant \frac{1}{2}\epsilon_{l_j}^{-1}.
\end{equation}

By $(\ref{pickj})$ and the choice of $j$, we have
\begin{equation}\label{pre}
\epsilon_{l_{j-1}}^2 > T(\tau)(\frac{\gamma}{(1+4l_{j-1}\ln\frac{1}{\epsilon_{l_{j-1}}} )^{\tau}})^{11}\epsilon_0^{'}(\frac{1}{l_{1}},\frac{1}{l_{2}}).
\end{equation}

Therefore, by $(\ref{cau})$ and $(\ref{pre})$, we get $\tilde{B}_1\in C^{k_1}(2\T^d,SL(2,\R))$ with
$$
\lVert \tilde{B}_1\rVert_{k_1}\leqslant \lVert B_1\rVert_{k_1} \lVert B_{l_j}\rVert_{k_1}\leqslant (1+\delta)\lVert B_{l_j}\rVert_{k_1}\leqslant \epsilon_{l_j}^{-1}\leqslant \gamma^{-20}C(\kappa,\tau,d,k,\lVert A\rVert)
$$
and
$$
\lvert \deg \tilde{B}_1\rvert = \lvert \deg B_{l_j}\rvert \leqslant 4l_{j}\ln\frac{1}{\epsilon_{l_j}}\leqslant \gamma^{-\frac{100}{k}}C(\kappa,\tau,d,k,\lVert A\rVert),
$$
which gives $(\ref{estb})$. For $(\ref{esta})$, one only need to notice that $\rho(\alpha, A_1)\neq0$, otherwise it will contradict to $\rho(\alpha, Ae^f)\in {\rm DC}_{\alpha}^d(\gamma,\tau)$. Thus $A_1$ can only be standard rotation in $SL(2,\R)$, which is the case of $(\ref{esta})$.

\textbf{(Rational case)} Again by $(\ref{sm3})$ and Proposition \ref{im}, there exists $B_{l_j}\in C^{\omega}_{\frac{1}{l_{j+1}}}(2\T^{d},SL(2,\R))$, $A_{l_j}\in SL(2,\R)$ and $\tilde{f}_{l_j}\in C^{k_0}(\T^{d},sl(2,\R))$ with $k_0= [\frac{k}{20}]$ such that
\begin{equation}
B_{l_j}(\theta+\alpha)(Ae^{f(\theta)})B^{-1}_{l_j}(\theta)=A_{l_j}e^{\tilde{f}_{l_j}(\theta)},
\end{equation}
with estimates
\begin{equation}\label{gi}
\lvert B_{l_j}(\theta)\rvert_{\frac{1}{l_{j+1}}}\leqslant \epsilon_{l_{j}}^{-\frac{2\sigma}{5}},\, \lvert \deg B_{l_j}\rvert\leqslant 4l_{j}\ln\frac{1}{\epsilon_{l_j}},\, \lVert \tilde{f}_{l_j}(\theta)\rVert_{k_0} \leqslant \epsilon_{l_j}^2.
\end{equation}

Since $\rho(\alpha, Ae^f)=\frac{\la m_0,\alpha\ra}{2}, m_0 \in \Z^d$ (we omit ``$\mod$" for simplicity), we have
$$
\rho(\alpha, A_{l_j}e^{\tilde{f}_{l_j}(\theta)})=\rho(\alpha, Ae^f)+\frac{\la \deg B_{l_j} ,\alpha\ra }{2}=\frac{\la m_0+\deg B_{l_j} ,\alpha\ra }{2}.
$$

Now, if we already have $m_0+\deg B_{l_j}=0$, i.e. $\rho(\alpha, A_{l_j}e^{\tilde{f}_{l_j}(\theta)})=0$, then by $(\ref{gi})$ and Lemma \ref{lem3.2}, there exists $B_2\in C^{k_1}(\T^d, SL(2,\R))$ with $k_1= [\frac{k_0}{20}]$ such that
\begin{equation}\label{main1}
B_2(\theta+\alpha)(A_{l_j}e^{\tilde{f}_{l_j}(\theta)})B_2^{-1}(\theta)=\bar{A},
\end{equation}
where $\bar{A}\in SL(2,\R)$ is a constant matrix and $B_2(\theta)$ is close to $Id$ (thus $\deg B_2=0$). Therefore we have
$$
\rho(\alpha,\bar{A})=\rho(\alpha, A_{l_j}e^{\tilde{f}_{l_j}(\theta)})+\deg B_2=0.
$$
Denote $\tilde{B}_2=B_2\circ B_{l_j}$, then $(\ref{estbb})$ follows.

On the other hand, if we have $m_0+\deg B_{l_j}\neq0$, then we can pick the smallest integer $j$ such that
\begin{flalign*}
\lvert \rho(\alpha,A_{l_j})\rvert \geqslant & \lvert \frac{\la m_0+\deg B_{l_j} ,\alpha\ra }{2}\rvert - \epsilon_{l_j}^2\\
\geqslant & \frac{\kappa}{2\lvert m_0+\deg B_{l_j}\rvert ^{\tau}}-\epsilon_{l_j}^2\\
\geqslant & \frac{\kappa}{4\lvert m_0+\deg B_{l_j}\rvert ^{\tau}}.
\end{flalign*}

Therefore, we can find a $P\in SL(2,\R)$, such that $\tilde{A}_{l_j}=PA_{l_j}P^{-1}$ is diagonal with estimates
$$
\lVert P\rVert \leqslant \sqrt{\frac{\lVert A_{l_j}\rVert}{\rho(\alpha,A_{l_j})}}\leqslant \sqrt{2\lVert A\rVert \times \frac{4\lvert m_0+\deg B_{l_j}\rvert ^{\tau}}{\kappa}}\leqslant \sqrt{8\kappa^{-1}\lVert A\rVert \lvert m_0+\deg B_{l_j}\rvert ^{\tau}}.
$$

Denote $\tilde{f}_{l_j}^{'}(\theta)=P\tilde{f}_{l_j}(\theta)P^{-1}$, we have
$$
\rho(\alpha, \tilde{A}_{l_j}e^{\tilde{f}_{l_j}^{'}(\theta)})=\rho(\alpha, A_{l_j}e^{\tilde{f}_{l_j}(\theta)})=\frac{\la m_0+\deg B_{l_j} ,\alpha\ra }{2}.
$$

Thus we can use $Q(\theta):=R_{-\frac{1}{2}\la m_0+\deg B_{l_j} ,\theta \ra }$ to conjugate $(\alpha, \tilde{A}_{l_j}e^{\tilde{f}_{l_j}^{'}(\theta)})$ to $(\alpha, \tilde{A}_{l_j}^{'}e^{\tilde{f}_{l_j}^{''}(\theta)})$:
$$
Q(\theta+\alpha)\tilde{A}_{l_j}e^{\tilde{f}_{l_j}^{'}(\theta)}Q^{-1}(\theta)=\tilde{A}_{l_j}^{'}e^{\tilde{f}_{l_j}^{''}(\theta)},
$$
with $\tilde{A}_{l_j}^{'}=Q(\theta+\alpha)\tilde{A}_{l_j}Q^{-1}(\theta) \in SL(2,\R)$ being still diagonal and $\lVert \tilde{A}_{l_j}^{'}\rVert \leqslant 2$.

Moreover, we have
$$
\rho((\alpha, \tilde{A}_{l_j}^{'}e^{\tilde{f}_{l_j}^{''}(\theta)}))=\rho(\alpha, \tilde{A}_{l_j}e^{\tilde{f}_{l_j}^{'}(\theta)})-\frac{1}{2}\la m_0+\deg B_{l_j} ,\theta \ra=0,
$$
and
\begin{flalign*}
\lVert \tilde{f}_{l_j}^{''}(\theta)\rVert_{k_0}& =\lVert Q(\theta)P\tilde{f}_{l_j}(\theta)P^{-1}Q^{-1}(\theta)\rVert_{k_0}\\
& \leqslant \lVert P\rVert^2 \lVert Q(\theta)\rVert_{k_0} \lVert\tilde{f}_{l_j}(\theta) \rVert_{k_0}\\
& \leqslant 8\kappa^{-1}\lVert A\rVert \lvert m_0+\deg B_{l_j}\rvert ^{\tau}\lvert \frac{1}{2}(m_0+\deg B_{l_j})\rvert^{k_0}\epsilon_{l_j}^2.
\end{flalign*}

Now by $(\ref{gi})$, we pick the smallest $j' \geqslant j$ such that
$$
8\kappa^{-1}\lVert A\rVert \lvert m_0+4l_{j'}\ln\frac{1}{\epsilon_{l_{j'}}}\rvert ^{\tau}\lvert \frac{1}{2}(m_0+4l_{j'}\ln\frac{1}{\epsilon_{l_{j'}}})\rvert^{k_0}\epsilon_{l_{j'}}^2\leqslant T(\tau)\kappa^{11}\epsilon_0^{'}(\frac{1}{l_{1}},\frac{1}{l_{2}})
$$
holds (such $j'$ exists since we choose $k\geqslant D_0\tau$ with $D_0$ large enough). Then by Lemma \ref{lem3.2} and by the same argument as above, $(\ref{estbb})$ is satisfied.
\end{pf}

\section{Spectral applications}
In this section, we concentrate on the $C^k$ quasi-periodic Schr\"{o}dinger operator $H_{V,\alpha,\theta}$. In other words, we will analyze this $C^k$ quasi-periodic Schr\"{o}dinger cocycle $(\alpha,S_E^{V})$:
$$
S_E^{V}(\theta)=\begin{pmatrix} E-V(\theta) & -1 \\ 1 & 0 \end{pmatrix}\in SL(2,\R),
$$
and use its reducibility results to obtain spectral properties of $H_{V,\alpha,\theta}$ and $\hat{L}_{V,\alpha,\varphi}$ (or $\hat{L}_{\lambda V,\alpha,\varphi}$ if we pick out $\lambda$ explicitly).
\subsection{Pure point of the dual operator}
In this subsection, we will prove Theorem $\ref{thm4.1}$. Before this, we first give a powerful criterion which was originally developed in \cite{AYZ} for the analytic case (and is a quantitative version of Theorem 4.1 of \cite{JK}), to establish a kind of equivalence between quantitative full measure reducibility of $C^k$ quasi-periodic Schr\"{o}dinger operator (cocycle) and pure point spectrum of its dual long-range operator.

\begin{Definition}
For any fixed $N\in\N,C>d,\varepsilon>0$, a normalized eigenfunction $u(\cdot)$ is said to be $(N,C,\varepsilon)$-good if
$$
\lvert u(m)\rvert \leqslant \varepsilon^{-1}\lvert m\rvert^{-C}, \ \  \lvert m\rvert\geqslant (1-\varepsilon)N.
$$
\end{Definition}
We label the $(N,C,\varepsilon)$-good eigenfunctions of $\hat{L}_{V,\alpha,\varphi}$ by $u^{\varphi}_j(\cdot)$ and denote the corresponding eigenvalue by $E^{\varphi}_j$. Moreover, denote
$$
\mathcal{E}^{\varphi}_{N,C,\varepsilon}=\{ E^{\varphi}_j\mid u^{\varphi}_j(\cdot)\,\, \mbox{is a $(N,C,\varepsilon)$-good normalized eigenfunction}\}
$$
and $\mathcal{E}_{C,\epsilon}(\varphi)=\bigcup_{N>0}\mathcal{E}^{\varphi}_{N,C,\varepsilon}$. Let $\mu^{pp}_{\delta_n,\varphi}$ be the spectral measure supported on $\mathcal{E}_{C,\epsilon}(\varphi)$ with respect to $\delta_n$, which means
$$
\langle\delta_n,\chi_{\mathcal{E}_{C,\epsilon}(\varphi)}(\hat{L}_{V,\alpha,\varphi})\delta_n\rangle=\int_{\sigma(\hat{L}_{V,\alpha,\varphi})}\chi_{\mathcal{E}_{C,\epsilon}(\varphi)}(E)d\mu^{pp}_{\delta_n,\varphi}(E).
$$

Now we can state our crucial criterion.
\begin{Proposition}\label{pro4.1}
Suppose that there exists $C$, such that for any $\delta>0$, there exists $\varepsilon>0$ and $N'$ sufficiently large such that for a.e. $\varphi\in \T$ and for all $N>N'$, we have
\begin{equation}\label{cru}
\#\{\mbox{linearly independent $(N,C,\varepsilon)$-good eigenfunctions}\}\geqslant (1-\delta)(2N)^d,
\end{equation}
then for a.e. $\varphi\in \T$, we have $\mu_\varphi=\mu_{\delta_n,\varphi}=\mu^{pp}_{\delta_n,\varphi}$ for any $n\in\Z^d$.
\end{Proposition}
\begin{pf}
Fix $\varphi\in \T$ such that $(\ref{cru})$ holds. Denote
$$
K^{\varphi}_{N',C,\varepsilon}=\{j\in\N \mid u^{\varphi}_j(\cdot)\,\, \mbox{is a $(N',C,\varepsilon)$-good eigenfunction} \}.
$$
$\#K^{\varphi}_{N',C,\varepsilon}$ is finite if $N',C,\varepsilon$ are fixed. For fixed $n\in \Z^d$, pick $N>\max\{N',\frac{\lvert 2n\rvert}{\varepsilon}\}$ and then we have
\begin{equation}\label{spe}
\sum_{\lvert m-n\rvert\leqslant N}\lvert u^{\varphi}_j(m)\rvert^2> \sum_{\lvert m\rvert\leqslant (1-\varepsilon)N}\lvert u^{\varphi}_j(m)\rvert^2>1-\frac{2^d}{\varepsilon(C-d+1)}N^{-C+d}
\end{equation}
for $(N,C,\epsilon)$-good eigenfuction $u^{\varphi}_j(\cdot)$.

Let $\tilde{\mu}^{pp}_{\delta_m,\varphi}=\tilde{\mu}^{pp}_{\delta_m,\varphi}(N,C,\varepsilon)$ be the truncated spectral measure supported on $\mathcal{E}^{\varphi}_{N,C,\varepsilon}$. Then by spectral theorem and $(\ref{spe})$, we have
\begin{flalign*}
\frac{1}{(2N)^d}\sum_{\lvert m-n\rvert \leqslant N}\lvert \mu^{pp}_{\delta_m,\varphi}\rvert > & \frac{1}{(2N)^d}\sum_{\lvert m-n\rvert \leqslant N}\lvert \tilde{\mu}^{pp}_{\delta_m,\varphi}\rvert\\
=& \frac{1}{(2N)^d}\sum_{\lvert m-n\rvert \leqslant N} \la P_{\mathcal{E}^{\varphi}_{N,C,\varepsilon}}\delta_m,\delta_m \ra\\
=& \frac{1}{(2N)^d}\sum_{\lvert m-n\rvert \leqslant N} \sum_{j\in K^{\varphi}_{N,C,\varepsilon}}\la P_{E_j^{\varphi}}\delta_m,\delta_m\ra\\
>& \frac{1}{(2N)^d}\sum_{\lvert m\rvert \leqslant (1-\varepsilon)N}\sum_{j\in K^{\varphi}_{N,C,\varepsilon}} \lvert u_j^{\varphi}(m)\rvert^2\\
>& \frac{1}{(2N)^d}\#K^{\varphi}_{N,C,\varepsilon}(1-\frac{2^d}{\varepsilon(C-d+1)}N^{-C+d})\\
>& (1-\delta)(1-\frac{2^d}{\varepsilon(C-d+1)}N^{-C+d}).
\end{flalign*}
Here $\lvert \mu^{pp}_{\delta_m,\varphi}\rvert :=\mu^{pp}_{\delta_m,\varphi}(\mathcal{E}_{C,\epsilon}(\varphi))$.

By the definition of $\hat{L}_{V,\alpha,\varphi}$ (\ref{1.2}), we obtain that the set $\mathcal{E}_{C,\epsilon}(\cdot)$ is the same for any $\phi'$ whose difference from $\phi$ is $\langle m,\alpha\rangle$, i.e. $\mathcal{E}_{C,\epsilon}(\varphi)=\mathcal{E}_{C,\epsilon}(\varphi+\la m,\alpha\ra)$. Then we have
$$
\frac{1}{(2N)^d}\sum_{\lvert m\rvert \leqslant N}\lvert \mu^{pp}_{\delta_n,\varphi+\la m,\alpha\ra}\rvert> (1-\delta)(1-\frac{2^d}{\varepsilon(C-d+1)}N^{-C+d}).
$$

Let $N$ goes to $\infty$, since $\delta$ is arbitrary small, then
$$
\int_{\T}\lvert \mu^{pp}_{\delta_n,\varphi}\rvert d\varphi=1,
$$
by the ergodic theorem of $\Z^d$ actions \cite{ps}. Thus for a.e. $\varphi \in \T$ and for any $n\in\Z^d$, $\mu_\varphi=\mu_{\delta_n,\varphi}=\mu^{pp}_{\delta_n,\varphi}$.
\end{pf}

Denote $\Theta_{\gamma}=\{\varphi \mid \varphi \in {\rm DC}_{\alpha}^d(\gamma,\tau) \}$. Note that $\lvert \bigcup_{\gamma>0}\Theta_{\gamma} \rvert=1 $ implies for any $\delta>0$, there exists $\epsilon_0$ such that if $\lvert \gamma\rvert\leqslant \epsilon_0$, then $\lvert \Theta_{\gamma}\rvert>1-\frac{\delta}{3}$. Again by ergodic theorem of $\Z^d$ actions \cite{ps}, we have
$$
\lim_{N_0\rightarrow \infty}\frac{1}{(2N_0)^d}\sum_{\lvert m\rvert\leqslant N_0}\chi_{\Theta_{\gamma}}(\varphi+\la m,\alpha \ra)=\int_{\T}\chi_{\Theta_{\gamma}}(\varphi)d\varphi.
$$

Thus for $N$ large enough if we take $N_0=N(1-\frac{\delta}{3})$, we have
\begin{equation}\label{4.3}
\#\{m \mid \varphi+\la m,\alpha\ra\in \Theta_{\gamma}, \lvert m\rvert \leqslant N(1-\frac{\delta}{3}) \}\geqslant (2N)^d(1-\delta).
\end{equation}

For any $\varphi \in \Theta_{\gamma}$, we choose $N'$ sufficiently large such that $(\ref{4.3})$ holds for $N>N'$.

As a direct application of our criterion, we prove Theorem \ref{thm4.1} in the following.
\begin{pf}\textbf{(Pure point spectrum part)}

We rewrite
$$
S_E^{V}(\theta)=\begin{pmatrix}E & -1 \\ 1 & 0 \end{pmatrix}+ \begin{pmatrix}-V(\theta) & 0 \\ 0 & 0 \end{pmatrix}.
$$

Under the assumptions of Theorem $\ref{thm4.1}$ and suppose that $\rho(\alpha, S_{E_m}^{\lambda V})=\varphi+\la m,\alpha \ra \in {\rm DC}_{\alpha}^d(\gamma,\tau)$, by Theorem $\ref{quanthm}$, then there exists $C_1=\gamma^{-20}C(\kappa,\tau,
d,k,$ $\lVert A\rVert)$, $C_2=\gamma^{-\frac{100}{k}}C(\kappa,\tau,d,k,\lVert A\rVert)$ and $B_m \in C^{k'}(2\T^d,SL(2,\R))$ with $k'=[\frac{k}{400}]$ such that
\begin{equation}\label{B}
B_m^{-1}(\theta+\alpha)S_{E_m}^{\lambda V}B_m(\theta)=R_{\varphi+\la m',\alpha \ra},
\end{equation}
with estimates:
\begin{equation}\label{C}
\lVert B_m\rVert_{k'}\leqslant C_1, \ \ \lvert m-m'\rvert \leqslant \lvert 2\deg B_m\rvert \leqslant 2C_2.
\end{equation}

Rewrite $(\ref{B})$ as
\begin{equation}
B_m^{-1}(\theta+\alpha)S_{E_m}^{\lambda V}B_m(\theta)= \begin{pmatrix} e^{2\pi i(\varphi+\la m',\alpha \ra)} & 0 \\ 0 & e^{-2\pi i(\varphi+\la m',\alpha \ra)} \end{pmatrix}.
\end{equation}

Also we write $B_m(\theta)=\begin{pmatrix} z_{11}(\theta) & z_{12}(\theta) \\ z_{21}(\theta) & z_{22}(\theta) \end{pmatrix}$, then direct computation shows
\begin{equation}\label{long}
(E_m-\lambda V(\theta))z_{11}(\theta)=e^{-2\pi i(\varphi+\la m',\alpha \ra)}z_{11}(\theta-\alpha)+e^{2\pi i(\varphi+\la m',\alpha \ra)}z_{11}(\theta+\alpha).
\end{equation}

Taking the Fourier transformation for $(\ref{long})$, we obtain
\begin{equation}
\sum_{k\in \Z^d}\lambda V_k {\hat{z} _{11}}(n-k)+2\cos 2\pi(\varphi+\la m',\alpha \ra+\la n,\alpha\ra)\hat{z}_{11}(n)=E_m \hat{z}_{11}(n),
\end{equation}
then $\hat{z}_{11}(n)$ is an eigenfunction since $z_{11}\in C^{k'}(2\T^d,\R)$.
To normalize $\hat{z}_{11}(n)$, we need the following inequality which was proved in \cite{AYZ}:
\begin{equation}\label{lemm}
\lVert \hat{z}_{11}\rVert_{l^2} \geqslant (2\lVert B_m \rVert_{C^0})^{-1}.
\end{equation}

Now normalizing $\hat{z}_{11}(n)$ by $u^{\varphi}_m(n)=\frac{\hat{z}_{11}(n+m')}{\lVert \hat{z}_{11}\rVert_{l^2}}$ and we will prove it is $(N,k',\varepsilon)$-good. Let
$$
\varepsilon \leqslant \min \{\frac{\delta}{6}-\frac{8C_2}{N}, \frac{1}{C_1^2}(\frac{\delta}{13})^{k'} \}.
$$

Since $u^{\varphi}_m(n)=u^{\varphi+\la m',\alpha\ra}_m(n-m')$, by $(\ref{C})$ and $(\ref{lemm})$, for $\lvert n\rvert \geqslant N(1-\varepsilon)$, we have
\begin{align*}
\lvert u^{\varphi}_m(n)\rvert=& \lvert u^{\varphi+\la m',\alpha\ra}_m(n- m')\rvert \\
\leqslant & \lVert B_m\rVert_{k'}^2\lvert n-m'\rvert^{-k'}\\
\leqslant & C_1^2 \lvert \lvert n \rvert - \lvert m\rvert-2C_2 \rvert^{-k'}\\
\leqslant & C_1^2 \lvert \frac{\frac{\delta}{3}-\varepsilon}{1-\varepsilon}\lvert n\rvert-2C_2 \rvert^{-k'}\\
\leqslant & C_1^2 \lvert \frac{\frac{\delta}{3}-\varepsilon}{2-2\varepsilon}\rvert^{-k'}\lvert n\rvert^{-k'}\\
\leqslant & \varepsilon^{-1}\lvert n\rvert^{-k'},
\end{align*}
which means $u^{\varphi}_m(\cdot)$ is $(N,k',\varepsilon)$-good and $k'>d$ is obvious since we choose $k>D_0\tau$ with $D_0$ large enough and $\tau>d-1$. By Proposition $\ref{pro4.1}$ and the estimates above, we conclude that for a.e. $\varphi\in \T$, $\hat{L}_{\lambda V,\alpha,\varphi}$ has pure point spectrum.

\textbf{(Polynomial SULE part)}

For any $\varphi \in \Theta_{\gamma}$, we have $\varphi+\la m,\alpha\ra\in \Theta_{\gamma(\lvert 2m\rvert+1)^{-\tau}}$, thus $\rho(\alpha, S_{E_m}^{\lambda V})=\varphi+\la m,\alpha \ra \in {\rm DC}_{\alpha}^d(\gamma(\lvert 2m\rvert+1)^{-\tau},\tau)$. By Theorem $\ref{quanthm}$ and same argument above, we have $C_1,C_2$ substituted by $\tilde{C}_1=(\gamma(\lvert 2m\rvert+1)^{-\tau})^{-20}C(\kappa,\tau,d,k,\lVert A\rVert)$ and $\tilde{C}_2=(\gamma(\lvert 2m\rvert+1)^{-\tau})^{-\frac{100}{k}}C(\kappa,\tau,d,k,\lVert A\rVert)$. Therefore, for any $n\in \Z^d$,
\begin{align*}
\lvert u^{\varphi}_m(n)\rvert=& \lvert u^{\varphi+\la m',\alpha\ra}_m(n- m')\rvert \\
\leqslant & \lVert B_m\rVert_{k'}^2\lvert n-m'\rvert^{-k'}\\
\leqslant & \gamma^{-40}(\lvert 2m\rvert+1)^{40\tau}C^2(\kappa,\tau,d,k,\lVert A\rVert)\lvert n-m'\rvert^{-k'}\\
\leqslant & \gamma^{-40}(\lvert 2m'\rvert+1)^{40\tau}(\lvert 2m-2m'\rvert+1)^{40\tau}C^2(\kappa,\tau,d,k,\lVert A\rVert)\lvert n-m'\rvert^{-k'}\\
\leqslant & \tilde{C}\lvert m'\rvert^{40\tau+\frac{4000\tau^2}{k}}\lvert n-m'\rvert^{-k'}.
\end{align*}

Now we pick the parameters $C=\tilde{C}$, $b=40\tau+{4000\tau^2}$ so that $b$ is independent of $k$ and $\tilde{k}=k'$ depending on $k$, then by definition \ref{defsemi}, $\hat{L}_{\lambda V,\alpha,\varphi}$ has polynomial semi-uniformly localized eigenvectors.

This finishes the proof of Theorem \ref{thm4.1}.
\end{pf}

\begin{Remark}
We need our quantitative full measure reducibility results to get more information about the eigenvectors $($``semi-uniformly localized"$)$. This is why we do not use Jitormiskaya and Kachkovskiy's method \cite{JK} to prove pure point spectrum.
\end{Remark}

\subsection{Generic Cantor spectrum}
In this subsection, we are dedicated to proving Thoerem \ref{THM4.2}.

To achieve this, we first give a brief review of Moser-P\"{o}schel argument in $C^k$ version, readers can refer to \cite{MP} and \cite{puig} for details in $C^{\omega}$ case. For definiteness, we put the proof of Proposition \ref{prop4.2} in the Appendix.

\begin{Proposition}\label{prop4.2}
Let $\alpha\in {\rm DC}_d(\kappa,\tau)$, $V \in C^k(\T^d, \R)$ with $k\geqslant D_0\tau$ for some large $D_0\in \N$. Assume that $(\alpha,S_E^{V})$ is $C^{k,k_1}$ reducible to $(\alpha, B)$ with $B=\begin{pmatrix} 1 & c \\ 0 & 1 \end{pmatrix}, c \in \R$. Let $W: \T^d\rightarrow \R$ be $C^{k}(\T^d, \R)$ and $t$ be real. $Z=\begin{pmatrix} z_{11} & z_{12} \\ z_{21} & z_{22} \end{pmatrix}: \T^d \rightarrow SL(2,\R)$ is the $C^{k_1}$ conjugation map.

If $(1) \,c\neq0 $ and $[Wz_{11}^2]\neq0,$

or $(2) \,c=0$ and $-[Wz_{11}z_{12}]^2+[Wz_{12}^2][Wz_{11}^2]<0$,

then $(\alpha, S_E^{V+tW})$ is uniformly hyperbolic provided $\lvert t\rvert>0$ is sufficiently small and $ct[Wz_{11}^2]<0$ if $(1)$ holds.
\end{Proposition}

As a straight corollary of Proposition $\ref{prop4.2}$ , we obtain
\begin{Corollary}
Same assumptions as Proposition \ref{prop4.2} and assume that $c\neq0$, then $E$ is at the endpoint of a noncollapsed spectral gap$($the right point if $c>0$ and the left one if $c<0$$)$.
\end{Corollary}
\begin{pf}
Take $W=1$, then $(\alpha, S_{E-t}^{V})$ is uniformly hyperbolic if $ct<0$ and $\lvert t\rvert$ is small enough. This means that there is an open spectral gap besides $E$($E$ is the right endpoint if $c>0$ and the left one if $c<0$).
\end{pf}
\begin{Remark}
For the case $c=0$, the corresponding energy lies at the endpoint of a collapsed gap as the derivative of the rotation number is nonzero. One can consult \cite{puig} for detailed proof.
\end{Remark}

The next corollary of Proposition $\ref{prop4.2}$ shows that when $c=0$, i.e. the Floquet matrix is the identity, one can ``open up" the collapsed gap by means of a generic $C^{k}$ perturbation.

\begin{Corollary}\label{cor4.2}
Under the same assumptions as Proposition \ref{prop4.2} and assume that $c=0$. If $W$ is a generic $C^{k}$ potential, then for $\lvert t\rvert>0$ sufficiently small, the spectrum $\sigma(\alpha,V+tW)$ has an open spectral gap with the label corresponding to $N(E,V,\alpha)$.
\end{Corollary}
\begin{pf}
By Proposition $\ref{prop4.2}$, for a perturbation $\tilde{W}$ satisfying
\begin{equation}\label{star}
-[\tilde{W}z_{11}z_{12}]^2+[\tilde{W}z_{12}^2][\tilde{W}z_{11}^2]<0,
\end{equation}
$(\alpha, S_E^{V+t \tilde{W}})$ is uniformly hyperbolic if $\lvert t\rvert>0$ is sufficiently small. This means $E$ lies in a spectral gap of $\sigma(\alpha, V+t\tilde{W})$. By the continuity of IDS we have
$$
N(E,V+ tW,\alpha)=N(E,V,\alpha)
$$
for $\lvert t\rvert$ small enough.
Denote
$$
y_1=\frac{1}{2}(z_{11}^2+z_{12}^2), \ \ y_2=\frac{1}{2}(z_{11}^2-z_{12}^2),\ \ y_3=z_{11}z_{12},
$$
then $(\ref{star})$ can be written as
$$
[\tilde{W}y_1]^2-[\tilde{W}y_2]^2-[\tilde{W}y_3]^2<0.
$$

We pick $E'$ such that
$$
[(E'+W)y_1]=E'[y_1]+[Wy_1]=0.
$$
$E'$ is uniquely defined as $[y_1]\neq0$. Then $E'+W$ satisfies $(\ref{star})$ if and only if
\begin{equation}\label{starr}
[(E'+W)y_1]^2-[(E'+W)y_2]^2-[(E'+W)y_3]^2<0.
\end{equation}

Since $[(E'+W)y_1]=0$, then $(\ref{starr})$ holds unless
$$
[(E'+W)y_2]=0 \ \ \mbox{and} \ \ [(E'+W)y_3]=0,
$$
which is equivalent to
$$
-[Wy_1][y_2]+[Wy_2][y_1]=0 \ \ \mbox{and} \ \ -[Wy_1][y_3]+[Wy_3][y_1]=0,
$$
which is a $C^{k}$ (dense and open, hence) generic condition.
Therefore, if $\lvert t\rvert>0$ is small enough, the spectrum $\sigma(\alpha,V+tW)$ has an open gap with
$$
N(E-tE',V+tW,\alpha)=N(E,V,\alpha).
$$
\end{pf}

Now we are able to prove Theorem $\ref{THM4.2}$.

\begin{pf} For Schr\"{o}dinger operator $H_{V,\alpha,\theta}$ with $V \in C^k(\T^d, \R)$, Gap Labeling Theorem \cite{BLT,MP} says that for any fixed spectral gap, there exists a unique $m\in \Z^d$ such that $2\rho(\alpha, S_E^V)={\la m,\alpha\ra} \mod \Z$ is locally constant in the gap. Besides, different $m$ corresponds to different gaps and one gap is labeled by one $m$. by Theorem $\ref{quanthm}$, for $E\in \Sigma_{V,\alpha}$ with $2\rho(\alpha, S_E^V)-{\la m,\alpha\ra} \in\Z$, $(\alpha, S_E^V)$ can be conjugated to the identity or a parabolic matrix (hyperbolic matrix is excluded since $E$ is in the spectrum). Moreover, since the rotation number is strictly monotonic in the spectrum, we obtain that there is a countable dense subset of energies in the spectrum such that $(\alpha, S_E^V)$ is $C^{k_1}$ ($k_1<k$ fixed) conjugated to the identity or a parabolic matrix. These energies lie at endpoints of gaps which can be collapsed (the identity case). However, by corollary $\ref{cor4.2}$, a generic and arbitrarily small $C^k$ perturbation opens the collapsed gaps easily. Since there is a countable number of gaps, then Theorem $\ref{THM4.2}$ follows.
\end{pf}

\section{Appendix}
\textbf{Proof of Proposition \ref{prop4.2}.}
By the assumptions, we have
$$
Z^{-1}(\theta+\alpha)S_E^{V}(\theta)Z(\theta)=B=\begin{pmatrix} 1 & c \\ 0 & 1 \end{pmatrix}.
$$
Direct computation gives
$$
Z^{-1}(\theta+\alpha)S_E^{V+tW}(\theta)Z(\theta)=B+tWP,
$$
with
$$
P=\begin{pmatrix} z_{11}z_{12}-cz_{11}^2 & -cz_{11}z_{12}+z_{12}^2 \\ -z_{11}^2 & -z_{11}z_{12} \end{pmatrix},
$$
here we denote $P=P(\theta)$ and $W=W(\theta)$ for simplicity.

By Proposition \ref{im}, $(\alpha, B+tWP)$ can be $C^{[\frac{k_1}{20}]}$ conjugated to $(\alpha, B+t[WP]+t^2 R)$ where $[\cdot]$ denote the averaging of a quasi-periodic function and $R$ depends $C^{[\frac{k_1}{20}]}$ on $t,\theta$.

Moreover, a direct computation shows that we can rewrite
$$
B+t[WP]+t^2 R= \exp(\tilde{B}_0+t\tilde{B}_1+t^2\tilde{R}),
$$
with
$$
\tilde{B}_0=\begin{pmatrix} 0 & c \\ 0 & 0 \end{pmatrix}, \ \ \tilde{B}_1=\begin{pmatrix} [Wz_{11}z_{12}]-\frac{c}{2}[Wz_{11}^2] & -c[Wz_{11}z_{12}]+[Wz_{12}^2] \\ -[Wz_{11}^2] & -[Wz_{11}z_{12}]+\frac{c}{2}[Wz_{11}^2] \end{pmatrix}
$$
and $\tilde{R}\in sl(2,\R)$ depending $C^{[\frac{k_1}{20}]}$ on $t,\theta$.

Let $D=\begin{pmatrix} d_1 & d_2 \\ d_3 & -d_1 \end{pmatrix}= \tilde{B}_0+t\tilde{B}_1$ and $d=det(D)=-d_1^2-d_2d_3$.

We first consider case $(1)\, c\neq0 $ and $[Wz_{11}^2]\neq0$. We can compute that $d=ct[Wz_{11}^2]+\mathcal{O}(t^2)$. Since $ct[Wz_{11}^2]<0$, then $d<0$. We construct
$$
Q=\begin{pmatrix} d_2 & d_2 \\ -d_1+\sqrt{-d} & -d_1-\sqrt{-d} \end{pmatrix},
$$
then $det(Q)=2c\sqrt{-ct[Wz_{11}^2]}+\mathcal{O}(t)$ and
$$
Q^{-1}DQ=\triangle=\begin{pmatrix} \sqrt{-d} & 0 \\ 0 & -\sqrt{-d} \end{pmatrix},
$$
thus $Q$ conjugates $(\alpha, \exp(\tilde{B}_0+t\tilde{B}_1+t^2\tilde{R}))$ to $(\alpha,\exp(\triangle+\tilde{S}_1))$ where
$$
\tilde{S}_1(t,\theta)=t^2 Q^{-1}\tilde{R}(t,\theta)Q, \ \ \tilde{S}_1 \sim \mathcal{O}(\lvert t\rvert^{\frac{3}{2}}), \forall \, \theta.
$$

Therefore,
$$
\triangle+\tilde{S}_1=\sqrt{-ct[Wz_{11}^2]}\big(\begin{pmatrix} 1 & 0 \\ 0 & -1 \end{pmatrix}+\mathcal{O}(\lvert t\rvert)\big),
$$
which implies $(\alpha, S_E^{V+tW})$ is uniformly hyperbolic provided $\lvert t\rvert>0$ is sufficiently small and $ct[Wz_{11}^2]<0$.

Now we turn to case $(2) \,c=0$ and $-[Wz_{11}z_{12}]^2+[Wz_{12}^2][Wz_{11}^2]<0$. Then
$$
D=t\begin{pmatrix} [Wz_{11}z_{12}] & [Wz_{12}^2] \\ -[Wz_{11}^2] & -[Wz_{11}z_{12}] \end{pmatrix}:=t\tilde{D},
$$
and $\tilde{d}=det(\tilde{D})=-[Wz_{11}z_{12}]^2+[Wz_{12}^2][Wz_{11}^2]<0$. Similarly as above, there exists a conjugation map $Q$, independent of $t$ and $\theta$, which transforms $(\alpha,\exp(t\tilde{D}+t^2 \tilde{R}))$ into $(\alpha,\exp(t\tilde{\triangle}+t^2 \tilde{S}_2))$ where
$$
\tilde{\triangle}=\begin{pmatrix} \sqrt{-\tilde{d}} & 0 \\ 0 & -\sqrt{-\tilde{d}} \end{pmatrix},\ \ \tilde{S}_2(t,\theta)=Q^{-1}\tilde{R}(t,\theta)Q.
$$

Therefore,
$$
t\tilde{\triangle}+t^2\tilde{S}_2=t\sqrt{-\tilde{d}}\big( \begin{pmatrix} 1 & 0 \\ 0 & -1 \end{pmatrix}+\frac{t}{\sqrt{-\tilde{d}}}\tilde{S}_2 \big),
$$
which implies $(\alpha, S_E^{V+tW})$ is uniformly hyperbolic provided $t\neq0$ is small enough. This finishes our proof.

\section{Acknowledgements}
We are indebted to Jiangong You and Qi Zhou for their enthusiastic help. This work is supported by ``the Fundamental Research Funds for the Central Universities'', \textquotedblleft Deng Feng Scholar Program B\textquotedblright of Nanjing University, Specially-appointed professor program of Jiangsu province, PTDC/MAT-PUR/29126/2017, Nankai Zhide Fundation and NSFC grant (11671192).


\begin{thebibliography}{99}
\bibitem{A0}A. Avila, \textit{Global theory of one-frequency Schr\'{o}dinger operators.} {\it Acta Math.} \textbf{215}, 1-54 (2015)
\bibitem{A1}A. Avila, \textit{Almost reducibility and absolute continuity I.} Preprint.
\bibitem{ABD}A. Avila; J. Bochi; D. Damanik, {\it Cantor spectrum for Schr\"{o}dinger operators with potentials arising from generalized skew-shifts.} {\it Duke Math. J.} \textbf{146}, 253-280 (2009)
\bibitem{AD}A. Avila; D. Damanik, \textit{Generic singular spectrum for ergodic Schr\"{o}-
dinger operators.} {\it Duke Math. J.} \textbf{130}, 393-400 (2005)
\bibitem{AFK}A. Avila; B. Fayad ; R. Krikorian, \textit{A KAM scheme for $SL(2,\R)$ cocycles with Liouvillean frequencies.} {\it Geom. Funct. Anal.} \textbf{21}, 1001-1019 (2011)
\bibitem{AJ}A. Avila; S. Jitomirskaya, \textit{The Ten Martini Problem.} {\it Ann. of Math.} \textbf{170}, 303-342 (2009)
\bibitem{AK}A. Avila; R. Krikorian, \textit{Reducibility or non-uniformly hyperbolicity for quasi-periodic Schr\"{o}dinger cocycles.} {\it Ann. of Math.} \textbf{164}, 911-940 (2006)
\bibitem{AYZ}A. Avila; J. You; Q. Zhou, \textit{Sharp phase transitions for the almost Mathieu operator.} {\it Duke Math. J.} \textbf{166(14)}, 2697-2718 (2017)
\bibitem{AYZ2}A. Avila; J. You; Q. Zhou, \textit{Dry ten Martini problem in the non-critical case.} Preprint.
\bibitem{AS}J. Avron; B. Simon, \textit{Almost periodic Schr\"{o}dinger operators. \uppercase\expandafter{\romannumeral2}: The integrated density of states.} {\it Duke Math. J.} \textbf{50}, 369-391 (1983)
\bibitem{BLT}J. Bellissard; R. Lima; D. Testard, \textit{Almost periodic Schr\"{o}dinger operators.} {\it Mathematics + physics. Vol.1 1-64 World Sci. Publishing, Singapore} (1985)
\bibitem{BG}J. Bourgain; M. Goldstein, \textit{On non-perturbative localization with quasi-periodic potential.} {\it Ann. of Math.} \textbf{152}, 835-879 (2000)
\bibitem{CCYZ}A. Cai; C. Chavaudret; J. You; Q. Zhou, \textit{Sharp H\"older continuity of the Lyapunov exponent of finitely differentiable quasi-periodic cocycles.} {\it Math. Z.} \textbf{291(3-4)}, 931-958 (2019)
\bibitem{chavaudret2}C. Chavaudret, \textit{Strong almost reducibility for analytic and gevrey quasi-periodic cocycles.} {\it Bull. Soc. Math. France} \textbf{141}, 47-106 (2011)
\bibitem{DJLS}R. del Rio; S. Jitomirskaya; Y. Last; B. Simon, \textit{Operators with singular continuous spectrum \uppercase\expandafter{\romannumeral4}. Hausdorff dimensions, rank-one perturbations, and localization.} {\it J. d' Analyse Math.} \textbf{69}, 153-200 (1996)
\bibitem{DS}E. Dinaburg; Ya. G. Sinai, \textit{The one-dimentional Schr\"{o}dinger equation with a quasi-periodic potential.} {\it Funct. Anal. Appl.} \textbf{9}, 279-289 (1975)
\bibitem{Eli92}L. H. Eliasson, \textit{Floquet solutions for the 1-dimensional quasi-periodic Schr\"{o}dinger equation.} {\it Commun. Math. Phys.} \textbf{146}, 447-482 (1992)
\bibitem{Eli97}L. H. Eliasson, \textit{Discrete one-dimensional quasi-periodic Schr\"{o}dinger operators with pure point spectrum.} {\it Acta Math.} \textbf{179}, 153-196 (1997)
\bibitem{FSW}J. Fr\"{o}hlich; T. Spencer; P. Wittwer, \textit{Localization for a class of one dimensional quasi-periodic Schr\"{o}dinger operators.} {\it Commun. Math. Phys.} \textbf{132}, 5-25 (1990)
\bibitem{GS}M. Goldstein; W. Schlag, \textit{On resonances and the formation of gaps in the spectrum of quasi-periodic Schr\"{o}dinger equations.} {\it Ann. of Math.} \textbf{173}, 337-475 (2011)
\bibitem{houyou}X. Hou; J. You, \textit{Almost reducibility and non-perturbative reducibility of quasiperiodic linear systems.} {\it Invent. Math.} \textbf{190}, 209-260 (2012)
\bibitem{Ji95}
S. Jitomirskaya,
\textit{Almost Everything About the Almost Mathieu Operator, \uppercase\expandafter{\romannumeral2}.}
Proceedings of \uppercase\expandafter{\romannumeral11} International Congress of Mathematical
Physics, Int. Press, 373-382 (1995)
\bibitem{J}S. Jitomirskaya, \textit{Metal-Insulator Transition for the almost Mathieu operator.} {\it Ann. of Math.} \textbf{150}, 1159-1175 (1999)
\bibitem{JK}S. Jitomirskaya; I. Kachkovskiy, \textit{$L^2$-reducibility and localization for quasiperiodic operators.} {\it Math. Res. Lett.} \textbf{23}, 431-444 (2016)
\bibitem{JL}S. Jitomirskaya; Y. Last, \textit{Anderson Localization for the Almost Mathieu Equation, \uppercase\expandafter{\romannumeral3}. Semi-Uniform Localization, Continuity of Gaps, and Measure of the Spectrum} {\it Commun. Math. Phys.} \textbf{195}, 1-14 (1998)
\bibitem{JLiu}S. Jitormiskya; W. Liu; \textit{Universal hierarchical structure of quasiperiodic eigenfuctions.} {\it Ann. of Math.} \textbf{187} no.3, 721-776 (2018).
\bibitem{Johnson}R. Johnson, \textit{Exponential dichotomy, rotation number, and linear differential operators with bounded coefficients.} {\it J. Differ. Equations} \textbf{61}, 54-78 (1986)
\bibitem{JM}R. Johnson; J. Moser, \textit{The rotation number for almost periodic potentials.} {\it Commun. Math. Phys.} \textbf{84}, 403-438 (1982)
\bibitem{LYZZ}M. Leguil; J. You; Z. Zhao; Q. Zhou, \textit{Asymptotics of spectral gaps of quasi-periodic Schr\"odinger operators.} arXiv:1712.04700.
\bibitem{MP}J. Moser; J. P\"{o}schel, \textit{An extension of a result by Dinaburg and Sinai on quasi-periodic potentials.} {\it Comment. Math. Helvetici.} \textbf{59}, 39-85 (1984)
\bibitem{ps}M. Pollicott; K. Schmidt, \textit{Ergodic theory of $\Z^d$ actions: proceedings of the Warwick symposium, 1993-4.} {\it Cambridge University Press} (1996).
\bibitem{puig}J. Puig, \textit{A nonperturbative Eliasson's reducibility theorem.} {\it Nonlinearity} \textbf{19}, 355-376 (2006)
\bibitem{Simon}B. Simon, \textit{Almost periodic Schr\"{o}dinger operators: A review.} {\it Adv. appl. Math.} \textbf{3}, 463-490 (1982)
\bibitem{Sinai}Ya. G. Sinai, \textit{Anderson localization for one-dimensional difference Schr\"{o}dinger operator with quasi-periodic potential.} {\it J. Stat. Phys.} \textbf{46}, 861-909 (1987)
\bibitem{WZhou}J. Wang; Q. Zhou, \textit{Reducibility results for quasiperiodic cocycles with Liouvillean frequency.} {\it J. Dynam. Differential Equations} \textbf{24}, 61-83 (2012)
\bibitem{WZ}Y. Wang; Z. Zhang, \textit{Cantor spectrum for a class of $C^2$ quasiperiodic Schr\"{o}dinger operators.} {\it Int. Math. Res. Noti.} \textbf{8}, 2300-2336 (2017)
\bibitem{youzhou}J. You; Q. Zhou, \textit{Embedding of analytic quasi-periodic cocycles into analytic quasi-periodic linear systems and its applications.} {\it Commun. Math. Phys.} \textbf{323}, 975-1005 (2013)
\bibitem{Young}L. Young, \textit{Lyapunov exponents for some quasi-periodic cocycles.} {\it Ergodic Theory Dynam. Systems} \textbf{17}, 483-504 (1997)
\bibitem{zehnder}E. Zehnder, \textit{Generalized implicit function theorems with applications to some small divisor problems:\uppercase\expandafter{\romannumeral1}.} {\it Comm. Pure Appl. Math.} \uppercase\expandafter{\romannumeral20}\uppercase\expandafter{\romannumeral8}, 91-140 (1975)
\end{thebibliography}
\end{document}